\documentclass[10pt]{amsart}
\usepackage{color}
\usepackage{verbatim}

\usepackage[square,compress,comma, numbers,sort]{natbib}
\usepackage[colorlinks=true, citecolor=green, linkcolor=blue]{hyperref}

\allowdisplaybreaks[4]

\def\H{\mathcal{H}}

\def\LT{\left}
\def\RT{\right}

\definecolor{c20}{rgb}{0.,0.7,0.}
\definecolor{c30}{rgb}{0.,0.,1.}
\definecolor{c40}{rgb}{1,0.1,0.7}
\definecolor{c50}{rgb}{1,0,0}
\definecolor{c60}{rgb}{1,0.9,0.1}

\def\b{\vk{b}}
\def\x{\vk{x}}

\def\SI{\Sigma}

\newcommand{\abs}[1]{\left\lvert #1 \right\rvert}

\newcommand{\pk}[1]{\mathbb{P} \left\{ #1 \right \} }

\newcommand{\R}{\mathbb{R}}

\newcommand{\N}{\mathbb{N}}
\newcommand{\inr}{\in \R}

\newcommand{\ldot}{,\ldots,}

\newcommand{\limit}[1]{\lim_{#1 \to   \infty}}

\newcommand{\BQN}{\begin{eqnarray}}
\newcommand{\EQN}{\end{eqnarray}}
\newcommand{\BQNY}{\begin{eqnarray*}}
\newcommand{\EQNY}{\end{eqnarray*}}

\newcommand{\BS}{\begin{sat}}
\newcommand{\ES}{\end{sat}}
\newcommand{\BT}{\begin{theo}}
\newcommand{\ET}{\end{theo}}
\newcommand{\BK}{\begin{korr}}
\newcommand{\EK}{\end{korr}}

\newcommand{\BD}{\begin{de}}
\newcommand{\ED}{\end{de}}

\newcommand{\BIT}{\begin{itemize}}
\newcommand{\EIT}{\end{itemize}}
\newcommand{\BDI}{\begin{description}}
\newcommand{\EDI}{\end{description}}

\newcommand{\BRM}{\begin{remarks}}
\newcommand{\ERM}{\end{remarks}}

\newcommand{\BEL}{\begin{lem}}
\newcommand{\EEL}{\end{lem}}

\newtheorem{theo}{Theorem}[section]
\newtheorem{sat}[theo]{Proposition}
\newtheorem{de}[theo]{Definition}
\newtheorem{lem}[theo]{Lemma}

\newtheorem{korr}[theo]{Corollary}

\newtheorem{remarks}[theo]{Remarks}

\newcommand{\nelem}[1]{{Lemma \ref{#1}}}

\newcommand{\netheo}[1]{{Theorem \ref{#1}}}

\newcommand{\COM}[1]{}

\newcommand{\QED}{\hfill $\Box$}

\newcommand{\kb}[1]{\boldsymbol{#1}}
\newcommand{\vk}[1]{\kb{#1}}

\def\IF{\infty}

\def\LT{\left}
\def\RT{\right}

\def\vn{\varepsilon}

\def\Del{\triangle}

\def\HAS{\mathcal{H}_I}

\def\K1#1{\textcolor{black}{#1}}

\def\EH#1{\textcolor{black}{#1}}

\def\ccj#1{\textcolor{black}{#1}}

\def\EHb#1{\textcolor{black}{#1}}

\def\K1#1{\textcolor{black}{#1}}

\def\k1#1{\textcolor{black}{#1}}

\def\e1#1{\textcolor{black}{#1}}
\def\kk#1{\textcolor{black}{#1}}

\def\LLc#1{\textcolor{black}{#1}}

\def\wHAS{\widetilde{\HAS}}

\def\ci{C_I}
\def\gt{\widehat g}
\def\ggt{ \widetilde{ g}}
\def\II{\mathbb I}
\def \clA{\mathcal A}

\begin{document}

\title{On the cumulative Parisian ruin of multi-dimensional Brownian motion risk models}

\author{Lanpeng Ji}
\address{Lanpeng Ji, School of Mathematics, University of Leeds, Woodhouse Lane, Leeds LS2 9JT, United Kingdom
}
\email{l.ji@leeds.ac.uk}

\bigskip

\date{\today}
 \maketitle

 {\bf Abstract:}
Consider a multi-dimensional Brownian motion which models \LLc{the surplus processes} of multiple lines of business of an insurance company. Our main result gives exact asymptotics for the cumulative Parisian ruin probability as the initial capital tends to infinity. An asymptotic distribution for the conditional cumulative Parisian ruin time is also derived. The obtained results on the cumulative Parisian ruin can be seen as generalizations of some of the results derived in \cite{DHJT18}. As a particular interesting case, the two-dimensional Brownian motion \LLc{risk} model is discussed in detail. 
%Our results suggest that  the company should not merge its two lines of business  if  the cumulative Parisian ruin probability is chosen as a  measure of risk. This provides an evidence in supporting the principle of  portfolio diversification.



 \section{Introduction}

 Consider an insurance company which operates simultaneously $d$ ($d\ge 1$) lines of business. It is  assumed that the surplus \LLc{processes} of these lines of business are described by a multi-dimensional risk %Brownian motion 
 model:
 \BQN\label{eq:U}
 \vk{U}(t)= \vk{u}+\vk{\mu}t -\vk{X}(t),\ \ \ t\ge 0,
 \EQN
 where $\vk{u}=(u_1,u_2,\ldots,u_d)^{\top}$, with $u_i\ge0$, is a (column) vector of initial capitals of these business lines, $\vk{\mu}=(\mu_1 \ldot \mu_d)^\top$, with $\mu_i>0$, is a  vector of net premium income rates, and 
$\vk{X}(t)=(X_1(t),X_2(t),\cdots,X_d(t))^\top, t\ge0$ is a vector of net loss processes by time $t$.

In recent years, there has been an \LLc{increasing} interest in risk theory in the study of multi-dimensional risk models with different stochastic processes modeling $\vk{X}(t), t\ge 0$; see, e.g., \cite{AsmAlb10} for an overview. 
%In the literature, there have been mainly two directions of investigation on this topic, namely, ruin probabilities and optimal dividend problem. We refer to \cite{ IB15,   SS16,  BP18, DHM18} and references therein for   ruin probabilities related studies, and \cite{ AM17, GSZ18, AM18} and references therein for studies on optimal dividend problem. %It turns out that 
In comparison with the well-understood 1-dimensional risk models, the study of multi-dimensional risk models is more challenging.

%%%%%%
\COM{
More recently,   two-dimensional Brownian motion models have drawn a lot of attention due to its tractability and practical relevancy. The optimal dividend problems for the two-dimensional Brownian motion models % with independent coordinates 
have been discussed in \cite{GS06,  G18, GSZ18}. The ruin probabilities  (which are actually exit problems) for multi-dimensional Brownian motion models have been discussed under different contexts; see \cite{LM07, DHM18, Met10, SW13, PuR08, Garbit14, MBor15, DHJT18} and the references therein. %Particularly,  it was shown in \cite{DMSW19} that the multi-dimensional Brownian model can serve as an approximation of a multi-dimensional classical risk model in a Markovian environment, making it an appealing model to work with.
}%%%%%%%%%5
We consider in this paper the multi-dimensional  Brownian  motion \LLc{risk
 model, i.e.,}
\BQN\label{eq:X}
\vk{X}(t)=A\vk{B}(t),\ \ \ \ t\ge 0,
\EQN
%denotes the approximated total claim amount process by time $t$. Here
where $A \in \R^{d\times d}$ is a non-singular matrix, and 
$\vk{B}(t)=(B_1(t),\ldots,B_d(t))^\top,t\ge0$ is a standard $d$-dimensional Brownian motion with independent coordinates. 
\LLc{Multi-dimensional Brownian motion risk models have drawn a lot of attention due to its tractability; see, e.g., \cite{DHM18,  GSZ18} and references therein}.

We shall investigate the {\it cumulative Parisian ruin} problem of the multi-dimensional Brownian motion risk model \eqref{eq:U} with \LLc{$\vk X$ defined by} \eqref{eq:X}. The cumulative Parisian ruin was first introduced by \cite{GR17} based on the occupation (or sojourn)  times of the 1-dimensional risk process. % and is due to its ties with the cumulative Parisian options; see, e.g., \cite{Hug99}. We refer to \cite{LkRe18} for some interesting discussions on the VaR-type risk measure derived from cumulative Parisian ruin. 
In the multi-dimensional setup the cumulative Parisian ruin time (at level $r>0$) is defined as
\BQN\label{eq:tauu1}
\tau_r(\vk{u}) :=\inf\LT\{t>0: \int_0^t \II{ (\vk{U}(s) <\vk 0) }ds>r    \RT\},
\EQN
where $\II{(\cdot)}$ is the indicator function, and the inequality for vectors $\vk{U}(s) <\vk 0$ is meant component-wise. 
As remarked in \cite{GR17} ``{\it the parameter $r$ could be interpreted as the length of a clock started at the beginning of the first excursion, paused when the process returns above zero, and resumed at the beginning of the next excursion, and so on.}".  Clearly, if $r$ is set to be $0$ one obtains the {\it simultaneous ruin time} $\tau_0(\vk{u})$ for the multi-dimensional Brownian motion risk model, i.e.,
\BQNY
\tau_0(\vk{u}):= \inf\LT\{t>0:  \vk{U}(t) <\vk 0     \RT\}=\inf\LT\{t>0:  U_i(t) < 0, \ \  \forall 1\le i\le d    \RT\},
\EQNY
%Additionally, the definition of cumulative Parisian ruin is due to its ties with the cumulative Parisian options; see, e.g., \cite{H99}.
which has been discussed  recently in  \cite{DHJT18} under a different context.

In this paper our primary focus is on the  infinite-time cumulative Parisian ruin probability defined as
\BQNY
\pk{\tau_r(\vk{u})<\IF}.
% \pk{\tau_r(u)<\IF}:==\mathbb{P}_{\vk{x}_u}\LT\{\int_0^\IF \II{ (\vk{X}(t)-  \boldsymbol{\mu}t\in \mathcal{U}) }dt>r \RT\}, %\ \ \ u\to\IF,
\EQNY
Note that in the 1-dimensional setup, we have from (5) in \cite{DLM18} (see also \cite{GR17})  that % for $\mu$ positive
\BQN\label{eq:Bm1}
\pk{\tau_r(u)<\IF}&=&\pk{\int_0^\IF\II{ (B_1(t)-  \mu_1 t>   u)} dt >r}\nonumber\\
&=& \LT(2(1+\mu_1^2 r) \Psi(\mu_1\sqrt r) -\frac{\mu_1\sqrt{2 r}}{\sqrt{\pi}} e^{-\frac{\mu_1 ^2 r}{2}}\RT)  e^{-2   \mu_1 u  } %, \ \ u\to \IF,
\EQN
 for all $u\in \R$, where $\Psi(s)$ is the standard normal survival function. %Apart from this 1-dimensional case, e

 It turns out that  explicit formula for the cumulative Parisian ruin probability in the multi-dimensional setup is  difficult to obtain. %and so far there have been no suitable tools  in deriving it. Therefore, 
In this case,  it is of  interest to derive some asymptotic results by \LLc{letting the initial capitals tend to infinity. % become large, by resorting to the extreme value theory; see, e.g., \cite{EKM97, Mik08, AsmAlb10}. %The asymptotic result is also of practical interest  
%As explained in \cite{Mik08} `` {\it ... the consideration of large initial capitals is not just a mathematical assumption but also an economic necessity, which is reinforced by the supervisory authorities. ...}".
%To this end, %in the rest of the paper make the following assumption: 
We  shall assume that}
 \BQNY
 \vk{u}=\vk \alpha u=(\alpha_1 u, \alpha_2 u, \ldots, \alpha_d u),\ \ \ \ \  \alpha_i>0, \ u\ge 0,
 \EQNY
 \LLc{and consider the  asymptotics of the cumulative Parisian ruin probability as $u\to \IF$.}
For simplicity, hereafter we denote
 \BQN\label{eq:tauu2}
 \tau_r(u):=\tau_r(\vk{u}),\ \ \ \ u\ge0.
 \EQN
 %%%%%%%%%%%%5
 \COM{We derive as our principle result the exact asymptotics of %as $u\to\IF$
 \BQN \label{eq:Pu}
 \pk{\tau_r(u)<\IF}:= \pk{\tau_r<\IF}= \mathbb{P}
\LT\{\int_0^\IF\II{(  \vk X(t)-\vk \mu t>\vk\alpha u)} dt >r \RT\},\ \ \ \ u\to\IF.
\EQN}
%%%%%%%%%%%%%55
Define %the following function 
\BQN \label{fgt}
g(t) = \frac{1}{t} \inf_{\vk{x} \ge \vk{\alpha}+\vk{\mu} t}  \vk{x}^\top \Sigma^{-1}  \vk{x}, \ \ t\ge 0,\ \quad \text{with}\ \ 
\Sigma=AA^\top,
\EQN
 where $1/0$ is understood as $\IF$.
%The main result  of \cite{DHJT18}
Our principal result presented in Theorem \ref{Thm1} shows that, for any $r\ge0,$
\BQN  \label{e.main2}
 \pk{\tau_r(u)<\IF}&=& \mathbb{P}
\LT\{\int_0^\IF\II{(  \vk X(t)-\vk \mu t>\vk\alpha u)} dt >r \RT\}\nonumber \\
&\sim & \ci\HAS(r) u^{\frac{1- m}{2}}e^{- \frac{ \inf_{t\ge 0} g(t)  }{2} u}, \ \ \ \ u\to\IF,
\EQN
where $\ci>0$, $m\in \N$  are known constants and $\HAS(r)$ is a counterpart of the
celebrated Pickands constant; explicit expressions of these constants will be displayed in Section \ref{Sec:main}.

\COM{ %%%%%%
Similarly as in \cite{DHJT18}, where the case $r=0$ is discussed,  we shall prove \eqref{e.main2} using the double-sum method combined with the theory of quadratic programming problem. 
%It is worth mentioning that when $y=0$ the above result in \eqref{e.main2} reduces to the main result derived in Theorem 3.1 of \cite{DHJT18}; particularly $\HAS(0)=\HAS$ the same as that defined therein. 
One of the difficulties for our proof comes from the fact that the  Bonferroni's inequality adopted in \cite{DHJT18} for the supremum functional cannot be used now, for which new inequalities for sojourn-type probabilities are proposed.  
 
With motivation from  \cite{GS06,  AAM17, AM18}, it is also of interest to see if merger of two lines of business will benefit the company. Our results for the two-dimensional Brownian motion models discussed in Section 4 suggest that when we consider the cumulative Parisian ruin probability as a  measure of risk,  it is better to keep operating two lines of business, as merger of two lines of business will  make the cumulative Parisian ruin probability larger. This provides an evidence in supporting the principle of  portfolio diversification.
 }
 %%%%%%%%%

As a by-product, we also derive in Theorem \ref{Thm1} %\ref{KorrRT}
 the asymptotic distribution of 
$$\tau_{r_2}(u) \lvert  \tau_{r_1}(u)<\IF, \ \ \ \ u\to \IF 
$$
for any $0\le r_1 \le r_2<\IF$.  
The approximation of the  above quantity is of interest in risk theory; it will provide us with some idea of when cumulative Parisian ruin actually occurred at level $r_2$  knowing that it has occurred at some level $r_1$. We refer to \cite{EKM97, AsmAlb10,  JR18} and references therein for related discussions on ruin times.

It is worth mentioning that there are some related interesting studies on the asymptotic properties of sojourn times above a high level of 1-dimensional (real-valued) stochastic processes; see, e.g., \cite{Ber74, Berman82, Berman92}. We refer to  \cite{DHPM17, DLM18}  for  recent developments. The multi-dimensional counterparts of this problem  are more challenging, and to the best knowledge of the author there has been no result in this direction. Our study on the cumulative Parisian ruin probability for the multi-dimensional Brownian motion risk models covers this gap in a sense by deriving some asymptotic properties of the sojourn times.  %is the first trial in this direction.
%Moreover, the obtained results extend the main results in \cite{DHJT18}. 
 
As an important illustration, the two-dimensional Brownian motion risk model is discussed in detail. Asymptotic results for the cumulative Parisian ruin probabilities and the conditional cumulative Parisian ruin times are obtained for the full range of the parameters involved in the model. %These results complement those obtained 
% Note that only a special two-dimensional Brownian motion risk model was discussed in \cite{DHJT18}.
 
The rest of this paper is organised as follows. In Section \ref{Prelim} we introduce some notation  and present some preliminaries, which are extracted from \cite{DHJT18}.
The main results  are presented in Section \ref{Sec:main}, followed by a discussion on the two-dimensional Brownian motion risk model in Section \ref{Sec:two-dim}. The technical proofs are displayed in  Section \ref{proofmain} and Section \ref{sec:Two-Dim}.

\section{Notation and Preliminaries}\label{Prelim}

We assume that all vectors are $d$-dimensional column vectors written in bold letters with  \K1{$d\ge2$}.  %Furthermore,
Operations with vectors are meant component-wise, e.g., $\lambda \x=\x \lambda = ( \lambda x_1 \ldot \lambda x_d)^\top$ for any $\lambda\inr, \x\inr^d$. 
%Recall that all vectors here are assumed to be $d$-dimensional column vectors written in bold letters with  \K1{$d\ge2$}.  
%Operations with vectors are meant component-wise. %, so % $\vk{x} \vk{y}=(x_1 y_1 \ldot x_d y_d)^\top$,
%$\abs{\vk{x}}=(\abs{x_1}   \ldot \abs{x_d })^\top$ and $\lambda \x=\x \lambda = ( \lambda x_1 \ldot \lambda x_d)^\top$ for any $\lambda\inr, \x\inr^d$. % In order to avoid confusion we shall write $\x \cdot \vk{y}$ instead of $\x \vk{y}$,
Further, we denote
$\vk{0}=(0\ldot 0)^\top \inr^d.$ %, \quad \vk{1}=(1\ldot 1)^\top \inr^d. $$
%and write $\vk{e}_j$ for the $m$-dimensional column vector with $j$th element 1 and all other elements equal 0.
%For any \EHb{non-empty subset} $\mathcal T\subset \R$, denote the inner set of $\mathcal T$ by $\mathcal T^o$ and its closure set by $\overline{\mathcal T}$.
If $I\subset \{1,\ldots,d\}$, then for a vector $\vk a\in\R^d$ we denote \EHb{by} $\vk{a}_I=(a_i, i\in I)$ a sub-block vector of $\vk a$. Similarly, if further $J \subset \{1\ldot d\}$, for  a matrix $M=(m_{ij})_{i,j\in \{1,\ldots,d\}}\in \R^{d\times d}$ we denote \EHb{by} $M_{IJ}\ccj{=M_{I,J}}=(m_{ij})_{i\in I, j\in J}$ \EHb{the} sub-block matrix of $M$
\EHb{determined by $I$ and $J$.} Moreover, write $M_{II}^{-1}=(M_{II})^{-1}$ for the inverse matrix of $M_{II}$ whenever it exists.

As we will see, the solution to the quadratic programming problem involved  in \eqref{fgt}  is the key to our discussions. We introduce
the next lemma stated in  \cite{HA2005} (see also  \cite{DHJT18}), which is important for several definitions in the sequel.

\BEL \label{AL}
Let
$M \in \R^{d \times d},d\ge 2$ be a positive definite  matrix. % with inverse $B=\Sigma^{-1}$.
If  $\vk{b}\in  \R ^d \setminus (-\infty, 0]^d $, then the quadratic programming problem
$$ P_M(\vk{b}): \text{Minimise $ \vk{x}^\top M^{-1} \vk{x} $ under the linear constraint } \vk{x} \ge \vk{b} $$
has a unique solution $\widetilde{\vk{b}}$ and there exists a unique non-empty
index set $I\subseteq \{1, \ldots, d\}$ such that
\begin{eqnarray}  \label{eq:IJi}
&&\widetilde{\vk{b}}_{I} = 
\vk{b}_{I}\not= \vk 0_I, \quad M_{II}^{-1} \vk{b}_{I}>\vk{0}_I,\\% >\vk{0}_V,\\
\text{and}&& \text{if}\  I^c =\{ 1, \ldots, d\} \setminus I \not=
\emptyset, \text{ then }
\widetilde {\vk{b}}_{I^c}%= - ((M^{-1} )_{I^cI^c})^{-1} (M^{-1})_{I^cI} \vk b_I
=  M_{I^cI}M_{II}^{-1} \vk{b}_{I}\ge \vk{b}_{I^c}.%,\\
\label{eq:hii}
%M_{VV}^{-1} \vk{b}_{I}>\vk{0}_I.
\end{eqnarray}
Furthermore,
\BQNY
%\label{eq:alfa} \max_{ \vk{w}\in [0,\infty)^d} \frac{(\vk{w}^\top \vk{b})^2 }
%{\vk{w}^\top M \vk{w}}&=&
\min_{\vk{x} \ge  \vk{b}}
\vk{x}^\top M^{-1}\vk{x} =  \widetilde{\vk{b}}^\top M^{-1} \widetilde {\vk{b}}   =
\vk{b}_{I}^\top M_{II}^{-1}\vk{b}_{I}>0. %, \\
% \label{eq:new}
%\vk{x}^\top M^{-1} \widetilde{ \vk{b}}&=& \vk{x}_I^\top M_{II}^{-1} \widetilde {\vk{b}}_I=
%\vk{x}_I^\top M_{II}^{-1}\vk{b}_I,  \quad
%\forall \vk{x}\in \R^d.
\EQNY
%If $\vk{b}= b \vk{1}, b\in (0,\infty)$, then $ 2 \le \sharp\{i: i\in I\}  \le d$.
\EEL
\begin{de}
The unique index set $I$ that defines the solution of the quadratic programming problem in question
will be referred to as the {\it essential index set}.
\end{de}

Consider the minimisation problem involved in  \eqref{fgt}, i.e., $\inf_{\vk{x} \ge \vk{\alpha}+\vk{\mu} t}  \vk{x}^\top \Sigma^{-1}  \vk{x}$. For any fixed $t\ge 0$, we define $\vk{b}(t)=\vk{\alpha}+\vk{\mu} t,$ and  let  $I(t)\subseteq \{1\ldot d\}$ be the  essential index set of the quadratic programming problem $P_{\Sigma}(\b(t))$.

Note that the two-layer minimisation problem in the exponent of  \eqref{e.main2}, i.e.,
\BQNY
\inf_{t\ge 0}  g(t) = \inf_{t\ge 0}\ \frac{1}{t} \ \inf_{\vk{x} \ge \vk{\alpha}+\vk{\mu} t}  \vk{x}^\top \Sigma^{-1}  \vk{x} 
\EQNY
has been solved in Lemma 2.2 in \cite{DHJT18}, with the aid of  \nelem{AL}.  More precisely, it is proved therein  the function $g(t), t\ge 0$ is convex and attains its unique minimum at some $t_0$. Let $I=I(t_0)$ be  the essential index set of the quadratic programming problem $P_\Sigma (\vk b)$ with
$\vk{b}=\b(t_0)= \vk \alpha + \vk \mu t_0.$
Then
\BQN\label{eq:t0}
t_{0}=   \sqrt{\frac{\vk{\alpha}_{I}^\top \Sigma^{-1}_{II}  \vk{\alpha}_{I}   }{\vk{\mu}_{I}^\top \Sigma^{-1}_{II}  \vk{\mu}_{I}  }
	}>0,
\EQN
and
\BQN \label{eq:intr1}
%\hat g:=
g(t_0) = \inf_{t\ge0}g(t) %\frac{1}{t} \inf_{\vk{v} \ge \vk{\alpha}+\vk{\mu} t}  \vk{v}^\top \Sigma^{-1}  \vk{v}
=\frac{1}{t_0}  \vk{b}^\top_{I} \Sigma_{II}^{-1}  \vk{b}_{I}.
\EQN
Hereafter, we shall use the notation $\vk{b}=\b(t_0)$, and use $I=I(t_0)$ for the {\it essential index set} of the quadratic programming problem $P_{\Sigma}(\b)$.  Furthermore,
let  $\widetilde{\vk{b}}$ be the \EHb{unique} solution
of $P_{\Sigma}(\b)$. %with the   \EH{essential index} set $I$.
If $I^c=\{1,\ldots,d\}\setminus I\neq\emptyset $, we define (cf. \eqref{eq:hii})
 {\it weakly essential index set} and   {\it unessential index set} by
\BQN\label{K.def}
K=\{j\in I^c: \widetilde{\vk{b}}_j=\Sigma_{jI}\Sigma_{II}^{-1}\vk{b}_I=\vk{b}_j\},
\quad
\text{  and } J=\{j\in I^c: \widetilde{\vk{b}}_j=\Sigma_{jI}\Sigma_{II}^{-1}\vk{b}_I>\vk{b}_j\}.
\EQN
%which will play certain roles in the asymptotic results.
%\cc{to be the weakly essential index set and inessential index set, respectively,}
As we shall see, the index set $I$ determines $m,$ $\inf_{t\ge0}g(t) $ and  $\HAS(r)$  in the asymptotics \eqref{e.main2},
whereas both  $I$ and $K$ determine the constant $\ci$. Moreover, the set $J$, whenever non-empty, contains indices that do
 not play any role in the asymptotic result, but it does appear in the proof (see \eqref{eq:EE}).

Next, define  for $t>0$
\BQNY
g_I(t):= \frac{1}{t}  \vk{b}(t)^\top_{I} \Sigma_{II}^{-1}  \vk{b}(t)_{I} =\frac{1}{t}    \vk{\alpha} ^\top_{I}  \Sigma^{-1}_{II}   \vk{\alpha} _{I} +
2 \vk{\alpha} ^\top_{I}  \Sigma^{-1}_{II}   \vk{\mu}  _{I}+
\vk{\mu}^\top_{I}  \Sigma^{-1}_{II}   \vk{\mu} _{I} t.
\EQNY
Clearly, by  %\nelem{AL} and 
\eqref{eq:intr1} we have 
$$\gt:=g(t_0)=g_I(t_0). $$
Furthermore, %we have by Taylor expansion
%\BQN\label{eq:gt0t1}
%g_I(t_0+t)= \gt  +\frac{\ggt}{2}t^2(1+o(1)),\ \ t\to 0,
%\EQN
%with
we denote $$
\ggt:=g_I''(t_0)=2 t_0^{-3} ( \vk{\alpha}_{I}^\top \Sigma^{-1}_{II}  \vk{\alpha}_{I}),
%\frac{( \vk{\alpha}_{I_j}^\top \Sigma^{-1}_{I_jI_j}  \vk{\alpha}_{I_j})^{3/2}} {(\vk{\mu}_{I_j}^\top \Sigma^{-1}_{I_jI_j}  \vk{\mu}_{I_j} )^{1/2}},
$$
which will appear in the definition of the constant $C_I$ in Section \ref{Sec:main}.
%For notational simplicity we shall set below
%\BQN\label{GTT}
% \gt= \inf_{t\ge 0} g(t)= g_I(t_0) , \quad  \ \  \ggt= g_I^{''}(t_0).
 %\EQN

\section{Main Results}\label{Sec:main}

We  introduce some constants that will appear in the main results. First we
write   
$$m:=\sharp\{i: i\in I\}\  \ge 1 $$
  for the number of elements of the essential index set  $I$. Further,  define the following constant (existence is confirmed in \netheo{Thm1})
\BQN\label{eq:HITT}
\HAS(r)  = \lim_{T\to\IF}\frac{1}{T}\HAS(r, T), 
\EQN
 \text{with} 
\BQN\label{eq:HIrT}
   \HAS(r,T)=\int_{\R^{m}}e^{\frac{1}{t_0}\vk{x}_I^\top\Sigma_{II}^{-1}\vk{b}_I}
\pk{\int_{t\in[0,T]}\II {((\vk{X}(t)-\vk{\mu}t)_I>\vk{x}_I)} dt>r }\,d\vk{x}_I,  \ \ r<T.
\EQN
Moreover, set
\BQNY
\ci:= \frac{1}{\sqrt{(2\pi t_0)^m \abs{\Sigma_{II}}}} \int_{\R} e^{-\ggt \frac{x^2}{4}} \psi(x)\,dx,
\EQNY
where $\abs{\Sigma_{II}}$ denotes the determinant of the matrix ${\Sigma_{II}}$, and for $x\inr$
\BQN\label{eq:psi}
\psi(x)=
\left\{\begin{array}{cc}
1, & \hbox{if } K=\emptyset,\\
\pk{\vk{Y}_K>\frac{1}{\sqrt {t_0}}(\vk{\mu}_K-\Sigma_{KI}\Sigma_{II}^{-1}\vk{\mu}_I)x}, & \hbox{if } K\neq\emptyset .
\end{array}\right.
\EQN
 Here the  index set $K$ is defined in (\ref{K.def}),  $\vk{Y}_{K}$ is  a Gaussian random vector with mean vector $\vk0_K$ and covariance matrix ${D}_{KK}$  given by 
$${D}_{KK}=\Sigma_{KK}-\Sigma_{KI}
\Sigma_{II}^{-1}\Sigma_{IK}.$$
%We also write 
%$$\vk{Y}_{K}\overset{d}\sim{\mathcal N}(\vk{0}_{K},{D}_{KK}).$$

The next theorem constitutes our main results. Its proof \K1{is} demonstrated in Section \ref{proofmain}.
\BT\label{Thm1} % \label{Thm1}
%Let $\vk \alpha, \vk \mu$ satisfy \eqref{alphamu} and let $\gt, \ggt$ be given by \eqref{GTT}. 
	Let $\tau_{r}(u)$ be defined  in  \eqref{eq:tauu2} (see also \eqref{eq:tauu1}). We have, for any $r\ge0,$ 
\BQN\label{eq:PuAsym}
 \pk{\tau_r(u)<\IF}  \sim  \ci\HAS(r)  u^{\frac{1- m}{2}}e^{- \frac{\gt}{2} u},\ \ \ \ u\to\IF,
\EQN
%as $u\to \IF$, 
where
\BQN\label{eq:HIbounds}
	  0<   \H_I(r) <\IF,\ \ \ \ \forall r\ge0.
\EQN
	% Let $\tau_{r}(u)$ be defined  in \eqref{eq:tauu1} and \eqref{eq:tauu2}, and let  the function $\psi$ be defined  in  \eqref{eq:psi}.   
	Moreover, we have, for any $0\le r_1\le r_2<\IF$ and any $s \in \R$,
	\BQN\label{eq:tau0u}
	\lim_{u\to\IF}\pk{\frac{\tau_{r_2}(u) -t_0 u} {\sqrt{ 2u /\ggt }} \le s  \Bigl \lvert \tau_{r_1}(u)<\IF}=\frac{ \H_I(r_2)   \int_{-\IF}^s e^{- \frac{x^2}{2}} \psi(\sqrt{2/\ggt }x)\,dx  }{  \H_I(r_1) \int_{-\IF}^\IF e^{-\frac{x^2}{2}} \psi(\sqrt{2/\ggt }x)\,dx }.
	\EQN  
	\ET 
\begin{remarks}
%We can see from the proof that the above two results are both valid if we relax the assumptions on the positivity of $\vk\alpha$ and $\vk mu$.

(a). If $d=1$, we have from \eqref{eq:PuAsym} that %Theorem \ref{Thm1} that 
\BQNY
\pk{\tau_r(u)<\IF}  \sim   \frac{1}{\mu_1}\H_{\{1\}}(r)   e^{- 2\alpha_1\mu_1 u},\ \ \ \ u\to\IF.
\EQNY
This together with \eqref{eq:Bm1} yields that
\BQN\label{eq:HBexpl}
\H_{\{1\}}(r) = \mu_1 \LT(2(1+\mu_1^2 r) \Psi(\mu_1\sqrt r) -\frac{\mu_1\sqrt{2 r}}{\sqrt{\pi}} e^{-\frac{\mu_1 ^2 r}{2}}\RT).
\EQN
%This has also been directly calculated in \cite{Ling18}.

\LLc{ (b). As in \cite{DHJT18} we can   check that %we are  interested in the case that
%\k1{$\limit{u}p_r(u)= 0$}, for which we shall assume that there exists some $1\le i\le d$ such that
the  results in %Theorems \ref{Thm1}  %are both valid under weaker conditions on  $\vk\alpha$ and $\vk \mu$.  More precisely, 
Theorem \ref{Thm1}  still hold for general $\vk\alpha, \vk \mu \in \R^d$  such that
%\BQNY % \label{alphamu}
	$\alpha_i>0,   \mu_i>0$ for some $1\le i\le d$.}
%\EQN

\end{remarks}
%%%%%%%%%%%%%%%%%%%%%%%
\section{Two-dimensional Brownian motion risk models}\label{Sec:two-dim}
\def\tzz{t_0^{(0)}}
\def\tzo{t_0^{(1)}}
\def\tzt{t_0^{(2)}}
In this section, we focus on the  two-dimensional Brownian motion risk models given by
\BQN \label{eq:Uu}
\vk U(t)=\vk \alpha v +\vk \mu t - A \vk B(t), \ \ t\ge 0,
\EQN 
with $\vk\alpha =(\alpha_1, \alpha_2)^\top>\vk 0$ and $\vk \mu =(\mu_1,\mu_2)^\top>\vk 0$ and 
%$$A =  
%\left(
%   \begin{array}{cc}
%     1 &0 \\
%     \rho & \sqrt{1-\rho^2 }\\
%   \end{array}
% \right), \  \quad \rho \in (-1,1).$$
 $$
  \Sigma= A A^{\top}=
\left(
   \begin{array}{cc}
     1 & \rho \\
     \rho & 1 \\
   \end{array}
 \right),\  \quad \rho \in (-1,1). 
$$
 
% Note that 
 We aim to find the asymptotics of  the cumulative Parisian ruin probability
%\BQNY
 %\pk{\tau_r(u)<\IF},% = \mathbb{P}\LT\{\int_0^\IF\II{(  \vk X(t)-\vk \mu t>\vk\alpha u)} dt >r \RT\}
 % \  \ u\to \IF, 
%\EQNY
and 
the asymptotic distribution of the conditional  cumulative Parisian ruin as $u$ tends to infinity,
%$$\tau_{r_2}(u) \lvert  \tau_{r_1}(u)<\IF, \  \ u\to \IF $$
for all the possible values of $\rho\in(-1,1)$, $\vk\alpha =(\alpha_1, \alpha_2)^\top>\vk 0$ and $\vk \mu =(\mu_1,\mu_2)^\top>\vk 0$.   
% in which we can observe  how different entries of the covariance matrix yield different scenarios of asymptotic behaviour.  

%From the main results presented in Theorem \ref{Thm3} below  we can observe  how different values of $\rho$ yield different scenarios of the asymptotic behaviour, which shows an interesting reduction of dimension phenomenon; see \cite{DHJT18} for more discussions on this phenomenon.
%Moreover, by comparing the cumulative Parisian ruin probabilities, we can conclude that   merger of two lines of business always make the risk larger, which does not benefit the   company.
 
%Without loss of generality we assume %consider the interesting case that, $\vk \mu=(\mu_1,\mu_2)>\vk 0,$ $\vk \alpha=(\alpha_1,\alpha_2)>\vk 0$ and

%%Then the considered sojourn time is defined as
%Applying the results in Theorem \ref{Thm1} w

In order to simplify the analysis,  we first do some  variable changes. Consider $\mu_1=ab$, with $a=1/\mu_1, b=\mu_1^2$. By the self-similarity of Brownian motion we derive that
\BQNY
 \tau_r(u) &=&
\inf\LT\{t\ge 0: \ \int_0^t \II
{\small \left( \begin{array}{l}
 X_1(s)-\mu_1 s >\alpha_1 u\\
  X_2(s)-\mu_2 s >\alpha_2 u
      \end{array}
      \right)}
ds >r \RT\} \\
&=&
\inf\LT\{t\ge 0: \ \int_0^t \II
{\small \left( \begin{array}{l}
 X_1(b^{-1}( bs))-  a (bs )>\alpha_1 u\\
  X_2(b^{-1}( bs))-\mu_2 b^{-1}( bs) >\alpha_2 u
      \end{array}
      \right)}
ds >r \RT\} \\
 &\overset{d}=&
\inf \LT\{t\ge 0: \  \int_0^ {\mu_1^2 t}\II
{\small \left( \begin{array}{l}
 X_1(s)-  s >\alpha_1\mu_1 u\\
  X_2(s)-\mu_2/\mu_1 s >\alpha_2\mu_1 u
      \end{array}
      \right)}
ds >\mu_1^2r \RT\},
\EQNY
where $\overset{d}=$ denotes equivalence in distribution.
%Then we have 
%$$ \vk X(t)=\LT(B_1(t), \ \rho B_1(t)+\sqrt{1-\rho^2} B_2(t)\RT)^\top,\ \ \ \ \ t\ge0.$$
Next we  denote
\BQN \label{eq:uvr}
v=\alpha_1\mu_1 u, \ \ \mu=\mu_2/\mu_1, \ \ \alpha=\alpha_2/\alpha_1, \ \ \widetilde{r} =\mu_1^2r,
\EQN 
and define
\BQNY
\tau_{\widetilde{r}}(v):=  \inf \LT\{t\ge 0:   \int_0^ { t}\II
{\small \left( \begin{array}{l}
 X_1(s)-  s >v\\
  X_2(s)-\mu  s >\alpha v
      \end{array}
      \right)}
ds > \widetilde{r} \RT\}.
\EQNY
Clearly, we have 
\BQN \label{eq:tauuv}
\tau_{\widetilde{r}}(v) \ \overset{d}= \ \mu_1^{2}  \cdot   \tau_r(u).
\EQN 
%%%%%%%
\COM{
\BQNY
\pk{\tau_r(u)<\IF} =
\mathbb{P}\LT\{\int_0^\IF\II
{\small \left( \begin{array}{l}
 X_1(t)-  t >v\\
  X_2(t)-\mu  t >\alpha  v
      \end{array}
      \right)}
dt >\widetilde{r} \RT\} =: \pk{\tau_{\widetilde{r}}(v)<\IF}.
\EQNY
Due to this equivalence, we have (with $\widetilde{r}_i =\mu_1^2r_i, i=1,2$)
\BQNY
\pk{\tau_r(u)<\IF}=\pk{\tau_{\widetilde{r}}(v)<\IF},\ \  \tau_{r_2}(u) \lvert  \tau_{r_1}(u)<\IF \ \overset{d}=\ \mu_1^{-2}  \cdot ( \tau_{\widetilde{r}_2}(v) \lvert  \tau_{\widetilde{r}_1}(v)<\IF).
\EQNY
} %%%%%%%%% 
Thanks to this equivalence, we can derive the results  for the cumulative Parisian ruin time $\tau_{ r }(u)$, %displayed in Theorem  \ref{Thm3}, 
by applying Theorem \ref{Thm1} to the cumulative Parisian ruin time $\tau_{\widetilde{r}}(v)$
of the auxiliary risk model  
\BQN\label{eq:Unew}
\widetilde{\vk U}(t)=\widetilde{\vk \alpha} v + \widetilde{\vk \mu} t - A \vk B(t), \ \ t\ge 0,
\EQN
 with (recall also \eqref{eq:uvr})
\BQNY
\widetilde{\vk\alpha} =(1, \alpha)^\top>\vk 0,\ \ \ \widetilde{\vk \mu} =(1,\mu)^\top>\vk 0. %,\ \ \ v=\alpha_1\mu_1 u. % \widetilde{r} =\mu_1^2r.
 %,\ \ \  \widetilde{r} =\mu_1^2r.
\EQNY
Note that  the auxiliary risk model   defined in \eqref{eq:Unew} is easier to analyse as it involves a smaller number of parameters (namely, $\rho, \alpha, \mu$) than the original risk model \eqref{eq:Uu}.

The main results of this section are displayed in \netheo{Thm3} below. From these results we can observe  how different values of $\rho$ yield different scenarios of the asymptotic behaviour, which shows an interesting reduction of dimension phenomenon; see also \cite{DHJT18} for discussions on this phenomenon. The proof of \netheo{Thm3} is deferred to Section \ref{sec:Two-Dim}.

%From the main results presented in Theorem \ref{Thm3} below  

% Recall the notation in \eqref{eq:uvr}. We have the following result.

\BT \label{Thm3} 
Consider the original two-dimensional Brownian motion models described in \eqref{eq:Uu}.
Recall the notation in \eqref{eq:uvr}. 
\begin{itemize}

\item[(i).] Suppose $\alpha$ and $\mu$ satisfy one of the following conditions:
\begin{itemize}
\item[(i.C1)] $\mu<1$ and $\alpha<1$,
\item[(i.C2)] $\mu<1$, $\alpha\ge 1$ and $\mu\le 1/\alpha$,
\item[(i.C3)]  $\mu\ge 1$, $\alpha<1$ and $\mu\le 1/\alpha$.
\end{itemize}
We have, for  any $r\ge0,$  
\begin{itemize}
\item[(i.R1).] If $-1<\rho < \frac{\alpha+\mu}{2}$, then, as $u\to\IF$
\BQN \label{eq:i.1-1}
\pk{\tau_r(u)<\IF} \sim \frac{\H_{\{1,2\}}(\mu_1^2 r)  }{\sqrt{\alpha_1\mu_1t_0^2\pi (1-\rho^2) \ggt}}   u^{-\frac{1}{2}}e^{- \frac{\gt}{2} \alpha_1\mu_1 u },\ \ \ \ 
\EQN 
where
\BQNY
&&t_0=\sqrt{\frac{1+\alpha^2-2\alpha  \rho }{1+\mu ^2-2\mu \rho}},\quad  \gt=\frac{2}{t_0}\frac{1+\alpha^2-2\alpha\rho}{1-\rho^2} +
\frac{2(1+\alpha\mu-\mu\rho-\alpha\rho)}{1-\rho^2},
%(\widetilde{\vk\alpha}^\top \Sigma^{-1}\widetilde{\vk\alpha} +\widetilde{\vk\alpha}^\top \Sigma^{-1}\widetilde{\vk\mu }t_0),
\quad  \ggt=2t_0 ^{-3}\frac{1+\alpha ^2-2\alpha \rho }{1-\rho^2},
\EQNY
and, for any $\widetilde r \ge 0$ 
$$
\H_{\{1,2\}}(\widetilde r) =\lim_{T\to\IF}\frac{1}{T} \int_{\R^{2}}e^{\LT(\frac{1-\rho\mu }{ 1-\rho^2 }+\frac{1=\rho\alpha}{t_0(1-\rho^2)}\RT)x_1 +\LT(\frac{\mu-\rho }{1-\rho^2}+\frac{\alpha-\rho}{t_0(1-\rho^2)}\RT)x_2 }
\pk{\int_{t\in[0,T]}\mathbb{I} 
{\small \left( \begin{array}{l}
 X_1(t)-  t >x_1\\
  X_2(t)-\mu  t >x_2
      \end{array}
      \right)}
dt>\widetilde r}\,d\vk{x}.
$$

Furthermore,   for any $0\le r_1\le r_2<\IF$ and any $s\in\R$
\BQN\label{eq:i.1-2}\lim_{u\to\IF}\pk{\frac{\tau_{r_2}(u) - t_0 \alpha_1/\mu_1 u} {\sqrt{ 2 \alpha_1 / (\mu_1^3\ggt)  u } } \le s  \Bigl \lvert \tau_{r_1}(u)<\IF }=\frac{\H_{\{1,2\}}(\mu_1^2 r_2)}{\H_{\{1,2\}}(\mu_1^2  r_1)}\Phi(s), 
\EQN % is approximated by a standard normal random variable.
where $\Phi(s)$ is the   standard normal distribution function.

\item[(i.R2).] If $ \rho=\frac{\alpha+\mu}{2 }$, then, as $u\to\IF$
\BQNY\label{eq:i.2-1}
 \pk{\tau_r(u)<\IF} \sim %\frac{\H_{\{2\}}(r) }{\sqrt{ 2\pi/\mu_2}}  \int_{\R}e^{-\frac{\mu_2^3}{2 }x^2}\Psi\LT(\frac{\mu_1-\rho\mu_2}{\sqrt{(1-\rho^2)/\mu_2}}x\RT)\, dx
\frac{1}{2 } \H_{\{1\}}(\mu_1^2r) \  e^{-2\alpha_1\mu_1 u},
\EQNY
where %the explicit expression for $\H_{\{1\}}(r)$ is available 
(cf. \eqref{eq:HBexpl})
$$
 \H_{\{1\}}(\widetilde r)=    2(1+ \widetilde r) \Psi( \sqrt {\widetilde r}) -\frac{ \sqrt{2 \widetilde r}}{\sqrt{\pi}} e^{-\frac{ \widetilde r}{2}}.
 %\lim_{T\to\IF}\frac{1}{T} \int_{\R}e^{2\mu_2 x} \pk{\int_{t\in[0,T]}\mathbb{I}_{((X_2(t)-  \mu_2 t) >x)}dt>r}\,dx.
$$
Furthermore, for any $0\le r_1\le r_2<\IF$ and any $s\in\R$
	\BQNY\label{eq:i.2-2}
\lim_{u\to\IF}	\pk{\frac{\tau_{r_2}(u) -\alpha_1/\mu_1u} {\sqrt{  \alpha_1/\mu_1^3 u   } } \le s  \Bigl \lvert \tau_{r_1}(u)<\IF}= \frac{\sqrt{2}\  \H_{\{1\}}(\mu_1^2r_2)}{\sqrt{\pi} \  \H_{\{1\}}(\mu_1^2 r_1)  }\ \int_{-\IF}^s e^{- \frac{x^2}{2}}  \Psi\LT(\frac{\mu -\rho }{\sqrt{(1-\rho^2) } }x\RT) \,dx,
	\EQNY
	where $\Psi(s)=1-\Phi(s)$ is the standard normal survival function.

\item[(i.R3).]  If $\frac{\alpha+\mu}{2}<\rho<1$, then, as $u\to\IF$
\BQNY\label{eq:i.3-1}
 \pk{\tau_r(u)<\IF} \sim  \H_{\{1\}}(\mu_1^2 r) \  e^{-2\alpha_1\mu_1 u},
\EQNY
and for any $0\le r_1\le r_2<\IF$ and any $s\in\R$
\BQNY\label{eq:i.3-2} 
\lim_{u\to\IF}\pk{\frac{\tau_{r_2}(u) -\alpha_1/\mu_1u} {\sqrt{  \alpha_1/\mu_1^3 u   } } \le s  \Bigl \lvert \tau_{r_1}(u)<\IF}= \frac{\H_{\{1\}}(\mu_1^2r_2)  }{ \H_{\{1\}}(\mu_1^2r_1)  } 
 \Phi(s).
 \EQNY

\end{itemize}

\item[(ii).]  Suppose $\alpha$ and $\mu$  satisfy one of the following conditions:
\begin{itemize}
\item[(ii.C1)] $\mu\ge 1$ and $\alpha\ge1$,
\item[(ii.C2)] $\mu<1$, $\alpha\ge 1$ and $\mu> 1/\alpha$,
\item[(ii.C3)]  $\mu\ge 1$, $\alpha<1$ and $\mu> 1/\alpha$.
\end{itemize}
We have, for  any $r\ge0,$  
\begin{itemize}
\item[(ii.R1).] If $-1<\rho < \frac{\alpha+\mu}{2\alpha\mu}$, then 
\eqref{eq:i.1-1} and 
\eqref{eq:i.1-2} hold.

\item[(ii.R2).] If $ \rho=\frac{\alpha+\mu}{2 \alpha \mu}$, then, as $u\to\IF$
\BQNY\label{eq:ii.2-1}
 \pk{\tau_r(u)<\IF} \sim %\frac{\H_{\{2\}}(r) }{\sqrt{ 2\pi/\mu_2}}  \int_{\R}e^{-\frac{\mu_2^3}{2 }x^2}\Psi\LT(\frac{\mu_1-\rho\mu_2}{\sqrt{(1-\rho^2)/\mu_2}}x\RT)\, dx
\frac{1}{2 \mu} \H_{\{2\}}(\mu_1^2 r) \  e^{-2\alpha_2\mu_2 u},
\EQNY
where %the explicit expression for $\H_{\{2\}}(r)$ is available 
(cf. \eqref{eq:HBexpl})
$$
\H_{\{2\}}(\widetilde r)= \mu\LT(2(1+ \mu^2\widetilde r) \Psi( \mu \sqrt {\widetilde r}) -\frac{ \mu\sqrt{2 \widetilde r}}{\sqrt{\pi}} e^{-\frac{\mu^2 \widetilde r}{2}} \RT).
$$
Furthermore, for any $0\le r_1\le r_2<\IF$ and any $s\in\R$
	\BQNY\label{eq:ii.2-2}
\lim_{u\to\IF}	\pk{\frac{\tau_{r_2}(u) -\alpha_2/\mu_2 u} { \sqrt{ \alpha_2/\mu_2^3 u   } } \le s  \Bigl \lvert \tau_{r_1}(u)<\IF}= \frac{\sqrt{2}\  \H_{\{2\}}(\mu_1^2r_2)}{\sqrt{\pi} \  \H_{\{2\}}(\mu_1^2r_1)  }\ \int_{-\IF}^s e^{- \frac{x^2}{2}}  \Psi\LT(\frac{\mu_1-\rho\mu_2}{\sqrt{(1-\rho^2) }\mu_2}x\RT) \,dx,
	\EQNY
	%where $\Psi(s)=1-\Phi(s)$ is the standard normal survival function. 

\item[(i.R3).]  If $\frac{\alpha+\mu}{2\alpha\mu}<\rho<1$, then, as $u\to\IF$
\BQNY\label{eq:ii.3-1}
 \pk{\tau_r(u)<\IF} \sim  \frac{1}{\mu} \H_{\{2\}}(\mu_1^2 r) \  e^{-2\alpha_2\mu_2 u},
\EQNY 
and for any $0\le r_1\le r_2<\IF$ and any $s\in\R$
\BQNY\label{eq:ii.3-2} 
\lim_{u\to\IF}	\pk{\frac{\tau_{r_2}(u) -\alpha_2/\mu_2 u} { \sqrt{ \alpha_2/\mu_2^3 u   } } \le s  \Bigl \lvert \tau_{r_1}(u)<\IF}= \frac{  \H_{\{2\}}(\mu_1^2r_2)}{  \H_{\{2\}}(\mu_1^2r_1)  }
 \Phi(s).
 \EQNY

\end{itemize}
\end{itemize}

\ET

\section{Proof\EHb{s} of Theorem \ref{Thm1} } \label{proofmain} %
In this section we  present the proof of Theorems \ref{Thm1}. We shall focus on the case where $r>0$, since the case where $r=0$ has been included in \cite{DHJT18}.

First, by the self-similarity of Brownian motion we have,  for any $u>0$,
%We are interested in the asymptotics of
\BQNY
 \pk{\tau_r(u)<\IF} = \pk{\int_0^\IF\II{(\vk{X}(t)-\vk{\mu} t>\vk{\alpha} u)} dt>r  } = \pk{ u \ \int_0^\IF\II{(\vk{X}(t)>\sqrt u (\vk{\alpha}+\vk{\mu} t))}dt >r}.
\EQNY
Next, we have  the following sandwich bounds
\BQN\label{eq:pr}
p_r(u)\le  \pk{\tau_r(u)<\IF} \le  p_r(u)+ r_0(u),
\EQN
where
\BQNY
p_r(u):=\pk{u\ \int_{t\in \Del_u}  \II{(  \vk{X}(t)>\sqrt u (\vk{\alpha}+\vk{\mu} t) )} dt>r    }, \ \ r_0(u):=\pk{u\ \int _{t \in \widetilde\Del_u } \II{(\vk{X}(t)> \sqrt u (\vk{\alpha}+\vk\mu t))} dt>0 },
\EQNY
with \EHb{(recall the definition of $t_0$ in \eqref{eq:t0})}
$$
\Del_u=\left[t_0-\frac{\ln(u)}{\sqrt{u}}, t_0+\frac{\ln(u)}{\sqrt{u}} \right],\quad
\widetilde\Del_u=\left[0,t_0-\frac{\ln (u)}{\sqrt u} \right]\cup \left[t_0+\frac{\ln (u)}{\sqrt u},\IF \right).$$

In order to convey the main ideas and to reduce the complexity, we shall prove the theorem %divide the %aforementioned proof into 
in several steps and finally we complete the proof  by putting all the arguments together.

\subsection{Analysis of $r_0(u)$}\label{chopping off}

This step is concerned with sharp upper bound for $r_0(u)$ when $u$ is large. Note that
$$
r_0(u)=\pk{ \exists _{t \in \widetilde\Del_u }  \vk{X}(t)> \sqrt u (\vk{\alpha}+\vk\mu t)  }.
$$
The following result is Lemma 4.1 in \cite{DHJT18} (there was a misprint with $\sqrt u$  missing, and in eq.(30) therein $u$ should be $\sqrt u$).

\BEL \label{PitI}
 For all large $u$  we have
%\footnote{EH: the term $Cu$ in the bounds above was redundant, since we
%	have an error of size $e^{ \varepsilon (\ln u)^2}$ which is much heavier. I deleted $u$} 
\BQN \label{eq:t0m}
r_0(u)\le \EHb{C} \sqrt u e^{- \frac{ \gt u }{2 } -  \LT(\frac{ \min(g{''}(t_0+), g{''}(t_0-)  )}{2}-\vn\RT) (\ln(u))^2}
\EQN\footnote{Note that in general $g{''}(t_0+)\neq  g{''}(t_0-)$; see Remark A.7 in \cite{DHJT18} for an example. }
holds for some constant  $ C>0$ \ccj{and some sufficiently small $ \vn>0$} which do not depend on $u$.
\EEL

\subsection{Analysis of $p_r(u)$}
%We   investigate   the  asymptotics of $p_r(u)$ as $u \to\IF$.
%\EH{In order  to} do so,  we use the method of Piterbarg to  partition the interval $\Del_u$.
Denote, for any fixed $T>0$ and $u>0$
\BQNY %\label{eq:DelNu}
\Del_{j;u}=\Del_{j;u}(T)= [t_0+j Tu^{-1}, t_0+(j+1) Tu^{-1}],\ \ \ -N_u\le j\le N_u,
\EQNY
where $N_u=\lceil T^{-1} \ln(u) \sqrt{u} \rceil$
\LLc{(here $\lceil x\rceil$} denotes the smallest integer larger than $x$).

Denote
$$
\mathcal A_{j,u}=u\ \int_{t\in\Delta_{j;u}} \II{(\vk{X}(t)>
\sqrt{u}(\vk{\alpha}+t\vk{\mu}))} dt,
$$
and define
\begin{equation*} %\label{eq:PLp}
p_{r,j;u}=\pk{ \mathcal A_{j,u} >r},\ \ 
p_{r,i,j;u}=\pk{\mathcal A_{i,u}>r,\ \mathcal A_{j,u}>r    }.
 \end{equation*}
It follows, using a similar idea as in \cite{DLM18}, that
\BQN\label{eq:thetaT}
p_r(u) &\le& \pk{ \sum_{j=-N_u}^{N_u} \clA_{j,u} >r} \nonumber\\
&=& \pk{ \sum_{j=-N_u}^{N_u} \clA_{j,u} >r, \ \  \text{there\ exists\ exactly\  one  } j \text{\ such\ that\ } \clA_{j,u} >0 }\nonumber\\
&&+\pk{ \sum_{j=-N_u}^{N_u} \clA_{j,u} >r, \ \  \text{there\ exist\ } i\neq j \text{\ such\ that\ } \clA_{i,u} >0\  \& \ \clA_{j,u} >0  }\nonumber\\
&\le & p_{1,r}(u)+\Pi_0(u),
\EQN
and by Bonferroni's inequality  
\BQN\label{eq:thetaT2}
p_r(u) &\ge& \pk{ \sum_{j=-N_u+1}^{N_u-1} \clA_{j,u} >r} \nonumber\\
 &\ge& \pk{ \exists -N_u+1\le j\le  N_u-1 \text{\ such \ that\ }\clA_{j,u} >r} \nonumber\\
 &\ge& p_{2,r}(u)-\Pi_0(u),
\EQN
where
\begin{equation*} %\label{eq:Theta}
p_{1,r}(u)= \sum_{j=-N_u}^{N_u}p_{r,j;u},\ \
p_{2,r}(u) = \sum_{j=-N_u+1}^{N_u-1}p_{r,j;u},\ \
\Pi_0(u) = \sum_{ -N_u\le i< j \le  N_u}p_{0,i,j;u}.
\end{equation*}

In the following two subsetions we shall focus on the analysis of $p_{i,r}(u), i=1,2$, and $\Pi_0(u)$, respectively.

\subsubsection{ \bf Analysis of the single sum $p_{i,r}(u), i=1,2$.}
%We shall focus on the asymptotics of $p_{1,r}(u)$, which will be easily seen to be asymptotically equivalent to $p_{2,r}(u)$ as $u\to\IF$.

First,  with the aid of Lemma 4.2 in \cite{DHJT18} we can check that %  it is not difficult to check that by the definition in \eqref{eq:HIrT}  
\BQN\label{eq:HHr0}
0<\H_I(r,T)\le \H_I(0,T)<\IF.
\EQN

%%%%%%%%%%%%%%%%%
\BEL\label{lem:PL} For any $T>0$ and  $r\in (0,T)$, we have as $u\to\infty$
\begin{eqnarray}\label{eq:PL0}
  p_{1,r}(u) \sim   p_{2,r}(u) \sim
\frac{1}{\sqrt{(2\pi t_0)^{m}|\Sigma_{II}|}  }  \frac{\H_I(r,T)}{T} u^{\frac{1- m}{2}}e^{-\frac{\gt u}{2}}\int_{\R} e^{-\frac{\ggt  x^2}{4}}\psi(x)\,dx,
\end{eqnarray}
 where $\psi(x)$ is given in \ccj{\eqref{eq:psi}.}

\EEL

\def \cu {c_{j;u}}

{\bf Proof:}  We shall focus on the asymptotics of $p_{1,r}(u)$, which is  easily seen to be asymptotically equivalent to $p_{2,r}(u)$ as $u\to\IF$.

Fix $T>0$. We shall prove the lemma in two steps.  In Step I we derive that \eqref{eq:PL0} holds for any $r\in (0,T)$ at which $\H_I(r,T)$,  as a function of $r$, is   continuous, and then in Step II we show that $\H_I(r,T), r\in (0,T)$  is actually continuous everywhere, implying that \eqref{eq:PL0} holds for all $r\in (0,T)$.

\underline{\it Step I}: The claim follows from similar arguments as in the proof of Lemma 4.3 in \cite{DHJT18}. 
%%%%%% 
% With the new form of \eqref{eq:Duu}, one could derive, using the same arguments as those of \cite{DHJT18}, that \eqref{eq:PL0} holds for any $r\in(0,T)$ at which $\H_I(r,T)$  is continuous; see also Theorem 5.1 of \cite{DHPM17} for related discussions. %Interested readers, DLM18} 
  %%%%%%%%%%%%%
By the independence and stationary increments property and the self-similarity of Brownian motion  we derive that
\begin{eqnarray*}
p_{r,j;u}&=&
\pk{u\int_{t\in[t_0+\frac{jT}{u}, t_0+\frac{(j+1)T}{u}]}
  \II{\LT(\vk{X}(t_0+\frac{jT}{u})+\vk{X}(t)-\vk{X}(t_0+\frac{jT}{u})>
\sqrt{u} (\vk{\alpha}+t\vk{\mu})\RT)} dt >r}
\\
&=&
\pk{\int_{t\in[0,T]} \II{\LT(\vk{Z}_{j;u}+\frac{1}{\sqrt{u}}({\vk{X}}(t)
-t\vk{\mu})>\sqrt{u}\vk{b}_{j;u}\RT)} dt>r  },
\end{eqnarray*}
 where
  $\vk{Z}_{j;u}$ is an independent of $\vk B$ Gaussian random vector with mean vector $\vk 0$ and covariance matrix  $\Sigma_{j;u}=\cu\Sigma$ with $\cu=\cu(T)=t_0+j T/u,$ and  
  $$\vk{b}_{j;u}=\vk{b}_{j;u}(T)=\b(t_0+\frac{jT}{u})=\b+\frac{jT}{u} \vk{\mu}.$$
Denote (recall $I^c =K\cup J$ in \eqref{K.def})
\begin{eqnarray*}\vk{Z}_K(t,\vk{x}_I)&=&(\vk{X}(t)-t\vk{\mu})_K-\Sigma_{KI}\Sigma_{II}^{-1}\vk{x}_I,\\
  \vk{Z}_J(t,\vk{x}_I)&=&(\vk{X}(t)-t\vk{\mu})_J-\Sigma_{JI}\Sigma_{II}^{-1}\vk{x}_I,
  \end{eqnarray*}
  and define  $\vk{Y}_{I^c }$ to be an independent of $\vk B$  Gaussian random vector with mean vector $\vk 0_{I^c}$ and covariance matrix $D_{I^cI^c}=\Sigma_{I^cI^c}-\Sigma_{I^cI} \Sigma_{II}^{-1}\Sigma_{II^c}$.
 % \BQNY
%\vk{Y}_{I^c } \  \overset{d}\sim \  \mathcal{N}(\vk 0_{I^c}, D_{I^cI^c}), \ \ \ D_{I^cI^c}=\Sigma_{I^cI^c}-\Sigma_{I^cI}
%\Sigma_{II}^{-1}\Sigma_{II^c},
%\EQNY
Using the same arguments as in  \cite{DHJT18} gives
\BQN \label{eq:rTu}
p_{1,r}(u)&=& \frac{u^{-m/2}}{\sqrt{(2\pi)^{m}\abs{\SI_{II} }}}
  \sum_{{-N_u}\le j\le N_u} \frac{1}{\cu^{m/2}}\exp\LT(-\frac{1}{2} u   g_I(t_0+\frac{jT}{u})  \RT)
 \int_{\R^{m}}  f_{j;u}(T,\vk x_I)P_{j;u}(r, T,\vk x_I)\,d\vk{x}_I\nonumber\\
 &=:&\frac{{1}}{T}\frac{1}{\sqrt{(2\pi)^{m}\abs{\SI_{II} }}} u^{(1-m)/2}  e^{ - \frac{ \hat g  u  }{2}  }
\EH{ R_{r,T}(u)},
\EQN
where
\BQNY
R_{r,T}(u)&=& \exp\LT( \frac{ \hat g  u  }{2}  \RT) \frac{T}{\sqrt{u}}\sum_{\kk{-N_u}\le j\le N_u}
\frac{1}{\cu^{m/2}}\exp\LT(-\frac{1}{2} u   g_I(t_0+\frac{jT}{u})  \RT)
\nonumber\\
&&\qquad\quad
\times \int_{\R^{m}}
f_{j;u}(T,\vk x_I)P_{j;u}(r, T,\vk x_I)\,d\vk{x}_I,
%\label{eq:bPL}
\EQNY
and
\BQNY
f_{j;u}(T,\vk x_I)&=&  \exp\LT(\frac{1}{\cu} \vk x_I^\top\Sigma_{II}^{-1} (\vk{b}_{j;u})_I  - \frac{1}{2u \cu} \vk{x}_I^\top\Sigma_{II}^{-1} \vk{x}_I    \RT) ,\label{eq:fju}\\
P_{j;u}(r, T,\vk x_I) &=&  \pk{
    \int_{t\in[0,T]} \II{(E(j, u, r, \vk x_I, t) )}
      dt},\label{eq:Pju}
\EQNY
with  the event $E(j, u, r, \vk x_I, t)$ defined as
\begin{eqnarray}
{E(j, u, r, \vk x_I, t)=\left(
      \begin{array}{l}
        (\vk{X}(t)-t\vk{\mu})_I>\vk{x}_I\\
        \sqrt{c_{j;u}} \vk Y_{K}
+\frac{1}{\sqrt{u}}\vk{Z}_K(t,\vk{x}_I) >\frac{jT}{\sqrt{u}}
(\vk{\mu}_K-\Sigma_{KI}\Sigma_{II}^{-1}\vk{\mu}_I)\\
\sqrt{c_{j;u}} \vk Y_{J}+\frac{1}{\sqrt{u}}
\vk{Z}_J(t,\vk{x}_I)
>\sqrt{u}(\vk{b}_J-\Sigma_{JI}\Sigma_{II}^{-1}\vk{b}_I+
(\vk{\mu}_J-\Sigma_{JI}\Sigma_{II}^{-1}\vk{\mu}_I)\frac{jT}{u})
      \end{array}
      \right).}   \label{eq:EE}
\end{eqnarray}
Using  similar argument as in Section 5.4  in  \cite{DHJT18} we can prove that
\BQNY  %\label{rT}
\limit{u} R_{r,T}(u)&=& t_0^{-m/2}\H_I(r, T)
 \int_{-\infty}^\infty e^{-\frac{g_I{''}(t_0) x^2}{4}}\psi(x)\,dx
\EQNY
holds for any $r\in(0,T)$ at which $\H_I(r,T)$  is continuous.  This together with \eqref{eq:rTu} yields that \eqref{eq:PL0} holds for any $r\in (0,T)$ at which $\H_I(r,T)$  is   continuous.

\underline {\it Step II}: We show that $\H_I(r,T)$  is  continuous at any point $r\in(0,T)$. Hereafter, let $r\in(0,T)$ be arbitrarily chosen and  fixed. We shall adopt an idea of \cite{DLM18}. %Note that $   \HAS(r,T)$ is right-continuous at $0$. It thus remains to prove the continuity of  $   \HAS(r,T)$  over $(0, T)$.  
 Recall
 \BQNY
   \HAS(r,T)=\int_{\R^{m}}e^{\frac{1}{t_0}\vk{x}_I^\top\Sigma_{II}^{-1}\vk{b}_I}
\pk{\int_{t\in[0,T]}\II {((\vk{X}(t)-\vk{\mu}t)_I>\vk{x}_I)} dt>r }\,d\vk{x}_I.
\EQNY
%With the aid of the dominated convergence theorem, the continuity of  $\HAS(r,T)$ at the fixed $r$ follows if 
We first show  that
  \BQN \label{eq:Pstar0}
\int_{\R^{m}}e^{\frac{1}{t_0}\vk{x}_I^\top\Sigma_{II}^{-1}\vk{b}_I}
\pk{\int_{t\in[0,T]}\II {((\vk{X}(t)-\vk{\mu}t)_I>\vk{x}_I)} dt=r }\,d\vk{x}_I = 0.
\EQN
To this end, we consider the probability space $(C_d([0,T]), \mathcal F, \mathbb P^*)$ which is induced by the multi-dimensional Brownian motion with drift $\{\vk{B}(t)-A^{-1}\vk \mu t, t\in[0,T]\}$, where $C_d([0,T])$ is the Banach space of all $d$-dimensional continuous vector functions over $[0,T]$, and $\mathcal F$ is the Borel $\sigma$-field of  $C_d([0,T])$. With the above notation, %it is sufficient to show that %for any  $r\in(0,T)$
 \eqref{eq:Pstar0} becomes
\BQN \label{eq:Pstar}
\int_{\R^{m}}e^{\frac{1}{t_0}\vk{x}_I^\top\Sigma_{II}^{-1}\vk{b}_I}
\mathbb{P}^*\LT\{\vk w \in C_d([0,T]):\ \int_{t\in[0,T]}\II {((A\vk{w}(t))_I>\vk{x}_I)} dt=r \right\}\,d\vk{x}_I = 0.
\EQN 
% where $\vk w \in C_d([0,T])$. 
 Denote, for the fixed $r,$ %any $r\in (0,T)$,
\BQNY
D_{\vk x_I}^{(r)}=\left\{\vk w \in C_d([0,T]: \int_{t\in[0,T]}\II {((A\vk{w}(t))_I>\vk{x}_I)} dt=r \right\},\ \ \vk x_I\in \R^m.
\EQNY
By the continuity of $\vk w$, one can see that
$$
D_{\vk x_I}^{(r)} \cap D_{\vk x'_I}^{(r)}=\emptyset,\ \ \ \vk x_I\neq \vk x'_I\in \R^m.
$$ 
Thus, for any finite number of distinct points $\vk x_I^{(1)}, \ldots, \vk x_I^{(N)}\in \R^m$ we have 
$$
\sum_{i=1}^N \mathbb P^*\{D^{(r)}_{{\vk x_I^{(i)}} }\} =\mathbb P^*\{  \cup_{i=1}^N  D^{(r)}_{{\vk x_I^{(i)}} }\}\le 1.
$$
This means that the set defined by
$$A^{(r)}_n=\{\vk x_I: \vk x_I\in \mathbb{R}^m \  \text{such that} \  \mathbb{P}^* \{D^{(r)}_{\vk x_I} \}   > {1}/{n}\}$$
consists of at most $n-1$ distinct $\vk x_I$'s. Therefore,
$$\{\vk x_I: \vk x_I\in \mathbb{R}^m \ \text{such that} \ \mathbb{P}^* \{D^{(r)}_{\vk x_I} \}   >0\}=\cup_{n=1}^\IF A^{(r)}_n$$
must be a countable set.  Consequently,
\BQNY
\int_{\R^{m}}e^{\frac{1}{t_0}\vk{x}_I^\top\Sigma_{II}^{-1}\vk{b}_I}
\mathbb{P}^*\LT\{ D^{(r)}_{\vk x_I}  \right\}\,d\vk{x}_I = 0,
\EQNY
which means that \eqref{eq:Pstar0} holds for the fixed $r\in (0,T)$. Next, for any small $\vn\in (0, r/2)$ we have 
$$\pk{\int_{t\in[0,T]}\II {((\vk{X}(t)-\vk{\mu}t)_I>\vk{x}_I)} dt>r\pm\vn }\le \pk{\int_{t\in[0,T]}\II {((\vk{X}(t)-\vk{\mu}t)_I>\vk{x}_I)} dt>r/2 }.
$$
Note that $\HAS(r/2,T)<\IF$; see \eqref{eq:HHr0}. Thus,  by the dominated convergence theorem we derive that
\BQNY
\lim_{\vn \downarrow 0} \HAS(r+\vn,T) =\HAS(r,T),
\EQNY
and 
\BQNY
\lim_{\vn \downarrow 0} \HAS(r-\vn,T) =\HAS(r,T) + \int_{\R^{m}}e^{\frac{1}{t_0}\vk{x}_I^\top\Sigma_{II}^{-1}\vk{b}_I}
\pk{\int_{t\in[0,T]}\II {((\vk{X}(t)-\vk{\mu}t)_I>\vk{x}_I)} dt=r }\,d\vk{x}_I,
\EQNY
which together with \eqref{eq:Pstar0} conclude the continuity of  $\HAS(r,T)$ at this $r$. Since such  $r$ was arbitrarily chosen in  $(0,T)$, we conclude that
$\H_I(r,T), r\in(0,T)$  is a continuous function.
 This completes the proof.
  \QED

 %%%%%%%%%%%%%%%%%%%%%%%%%%%%%%%%%%%%%%%%%%%%%%%%%%%%%%%%%%%%%5

\def\wHAS{\widetilde{\HAS}}

\subsubsection{\bf Estimation of the double-sum $\Pi_0(u)$.} %\label{DSM}
In this subsection we shall focus on asymptotic upper bounds of % that   and then $T\to\IF$
$\Pi_0(u)$, as $u\to\IF$.
%for large $u,T.$
Note that 
\BQN \label{eq:doubsum1}
\Pi_0(u)=\sum_{ -N_u\le i< j \le  N_u}p_{0,i,j;u} %=\sum_{ -N_u\le i< j \le  N_u}\pk{\exists_{t\in \Del_{i;u}}  \vk{X}(t)>\sqrt u (\vk{\alpha}+\vk{\mu} t),  \ \exists_{t\in \Del_{j;u}} \vk{X}(t)>\sqrt u (\vk{\alpha}+\vk{\mu} t)    }.
=\underset{j=i+1}{\sum_{-N_u\le i<j\le N_u}p_{0,i,j;u}}+
\underset{j>i+1}{\sum_{-N_u\le i<j\le N_u}p_{0,i,j;u}}=:\Pi_{0,1}(u)+\Pi_{0,2}(u).
\EQN 
Since
\BQNY
p_{0,i,j;u}  = \pk{\exists_{t\in \Del_{i;u}}  \vk{X}(t)>\sqrt u (\vk{\alpha}+\vk{\mu} t),  \ \exists_{t\in \Del_{j;u}} \vk{X}(t)>\sqrt u (\vk{\alpha}+\vk{\mu} t)    },
\EQNY
we obtain from (52) in \cite{DHJT18} that
\BQN
  \lim_{u\to\IF}\frac{\Pi_{0,1}(u)}{u^{(1-m)/2} \exp\LT(-\frac{ \gt   }{2} u\RT) } 
= Q_1  \LT(\frac{2\HAS(0,T)}{T}- \frac{\HAS(0,2 T)}{T}\RT) 
\EQN
for some  constant $Q_1>0$ which does not depend on $T$.
Similarly, %as in  \cite{DHJT18} we have
\BQN\label{eq:doubsum2}
 \lim_{u\to\IF}\frac{\Pi_{0,2}(u)}{u^{(1-m)/2} \exp\LT(-\frac{ \gt   }{2} u\RT) }\le Q_2 T \sum_{j\ge 1} \exp\Bigl(- \frac{\gt   }{8t_0} (jT)\Bigr) 
\EQN
holds with some constant $Q_2>0$   which does not depend on $T$.

Now we are ready to present the proof of \netheo{Thm1}. 

\underline{\bf Proof of \eqref{eq:PuAsym} and \eqref{eq:HIbounds}.}
We have from  \eqref{eq:pr}-%, \eqref{eq:t0m}, \eqref{eq:thetaT}, \eqref{eq:thetaT2}, 
\eqref{eq:PL0} and \eqref{eq:doubsum1}-\eqref{eq:doubsum2} that, for any $T_1,T_2>0$
\BQN
\limsup_{u\to\IF}\frac{ \pk{\tau_r(u)<\IF} }{ \ci u^{\frac{1- m}{2}}e^{- \frac{\gt}{2} u}} &\le&   \frac{\H_I(r,T_1)}{T_1}  +  Q_1  \LT(\frac{2\HAS(0,T_1)}{T_1}- \frac{\HAS(0,2 T_1)}{T_1}\RT) + Q_2 T_1 \sum_{j\ge 1} \exp\Bigl(- \frac{\gt   }{8t_0} (jT_1)\Bigr), \label{eq:up}\\
\liminf_{u\to\IF}\frac{ \pk{\tau_r(u)<\IF} }{ \ci u^{\frac{1- m}{2}}e^{- \frac{\gt}{2} u}} &\ge&   \frac{\H_I(r,T_2)}{T_2}  -  Q_1  \LT(\frac{2\HAS(0,T_2)}{T_2}- \frac{\HAS(0,2 T_2)}{T_2}\RT) - Q_2 T_2 \sum_{j\ge 1} \exp\Bigl(- \frac{\gt   }{8t_0} (jT_2)\Bigr).\label{eq:low}
\EQN
Note that it has been shown in \cite{DHJT18} that 
$$
\H_I(0)= \lim_{T\to\IF}  \frac{\H_I(0,T)}{T} <\IF.
$$
Letting $T_2\to\IF$ in \eqref{eq:low}, with $T_1$ in \eqref{eq:up} fixed, we have
\BQNY
 \limsup_{T\to\IF}  \frac{\H_I(r,T)}{T} <\IF.
\EQNY
Furthermore, letting $T_1\to\IF$ we conclude that
\BQNY
\liminf_{T\to\IF}  \frac{\H_I(r,T)}{T} =  \limsup_{T\to\IF}  \frac{\H_I(r,T)}{T} <\IF.
\EQNY
Therefore, it remains to prove that
\BQN\label{eq:liminfH}
\liminf_{T\to\IF}  \frac{\H_I(r,T)}{T} >0
\EQN
holds. %Next we prove \eqref{eq:liminfH}. 
To this end, first note that
\BQNY
 \pk{\tau_r(u)<\IF} &\ge& p_r(u)\ge  \pk{ \sum_{j=-N_u+1; j\in\{2k:k\in \mathbb Z\}}^{N_u-1} \clA_{j,u} >r} \nonumber\\
 &\ge& \pk{ \exists -N_u+1\le j\le  N_u-1,  j\in\{2k:k\in \mathbb Z\} \text{\ such \ that\ }\clA_{j,u} >r} \nonumber\\
 &\ge& p_{3,r}(u)-\widetilde \Pi (u),
\EQNY
where 
\begin{equation*} 
p_{3,r}(u) = \sum_{j=-N_u+1; j\in\{2k:k\in \mathbb Z\}}^{N_u-1}p_{r,j;u},\ \ \ \ 
\widetilde\Pi(u) = \sum_{ -N_u\le i< j \le  N_u; i,j\in\{2k:k\in \mathbb Z\}}p_{0,i,j;u}.
\end{equation*}
Similar augments as in the derivation of \eqref{eq:low} gives that, for some $T_3>0,$
\BQNY
\liminf_{u\to\IF}\frac{ \pk{\tau_r(u)<\IF} }{ \ci u^{\frac{1- m}{2}}e^{- \frac{\gt}{2} u}} &\ge&   \frac{\H_I(r,T_3)}{2 T_3}  - Q_3 T_3 \sum_{j\ge 1} \exp\Bigl(- \frac{\gt   }{8t_0} (jT_3)\Bigr)
\EQNY
holds with some constant $Q_3>0$  which does not dependent on $T_3$. This together with \eqref{eq:up} yields that
\BQNY
\liminf_{T_1\to\IF}  \frac{\H_I(r,T_1)}{T_1} %\ge \liminf_{u\to\IF}\frac{ \pk{\tau_r(u)<\IF} }{ \ci u^{\frac{1- m}{2}}e^{- \frac{\gt}{2} u}} 
&\ge&   \frac{\H_I(r,T_3)}{2 T_3}  - Q_3 T_3 \sum_{j\ge 1} \exp\Bigl(- \frac{\gt   }{8t_0} (jT_3)\Bigr)\\
&\ge&   \frac{\H_I(r, r+1)}{2 T_3}  - Q_3 T_3 \sum_{j\ge 1} \exp\Bigl(- \frac{\gt   }{8t_0} (jT_3)\Bigr),
\EQNY
holds for all $T_3 \ge r+1,$ where the last inequality follows since $\H_I(r,T)$ as a function of $T$ is non-decreasing. Since for sufficiently large $T_3$ the right-hand side of the above formula is positive, we conclude that \eqref{eq:liminfH} is valid. 
Thus, the proof of \eqref{eq:PuAsym} and \eqref{eq:HIbounds} is complete. 
 \QED

\def\cu{c_{\lambda,u}}

\underline{\bf Proof of \eqref{eq:tau0u}.}
We have, for any $s\in\R$
\begin{eqnarray}\label{eq:taur1r2}
\pk{\frac{\tau_{r_2}(u) -t_0 u}{\sqrt u} \le s  \big\lvert \tau_{r_1}(u)<\IF}&=&
\frac
{\pk{\frac{\tau_{r_2}(u) -t_0 u}{\sqrt u} \le s ,\tau_{r_1}(u)<\IF}}
{\pk{\tau_{r_1}(u)<\infty}}\nonumber
\\
&=&\frac{\pk{ {\tau}_{r_2}(u)\le ut_0+\sqrt{u}s}}{\pk{\tau_{r_1}(u)<\infty}}\nonumber\\
&=&\frac{\pk{u \ \int_0^{ t_0  +s/ \sqrt u} \II{(\vk X(t)>(\vk \alpha +\vk \mu t)\sqrt u)} dt>r_2}}{\pk{\tau_{r_1}(u)<\infty}} .
\end{eqnarray}
Furthermore, using the same arguments as in the proof of \eqref{eq:PuAsym} %and Theorem 3.3 in \cite{DHJT18}, 
we can show,  as $u\to\IF$
\BQNY
&&\pk{u \ \int_0^{ t_0  +s/ \sqrt u} \II{(\vk X(t)>(\vk \alpha +\vk \mu t)\sqrt u)} dt>r_2}\\
&&\sim
\pk{u \ \int_{t_0-\ln(u)/\sqrt u}^{ t_0  +s/ \sqrt u} \II{(\vk X(t)>(\vk \alpha +\vk \mu t)\sqrt u)} dt>r_2} \\ 
&&\sim
\frac{\HAS(r_2)}{\sqrt{(2\pi t_0)^m \abs{\Sigma_{II}}}} \int_{-\IF}^s e^{-\ggt \frac{x^2}{4}} \psi(x)\,dx \,  u^{\EH{\frac{1- m}{2}}}e^{- \frac{\gt}{2} u}.
\EQNY
Consequently, by plugging the above asymptotics and \eqref{eq:PuAsym} into \eqref{eq:taur1r2} and rearranging, we obtain \eqref{eq:tau0u}.  Thus, the proof is complete.  
\QED

\section{Proof of Theorem \ref{Thm3}}\label{sec:Two-Dim}

%Due to this equivalence, we can derive the results  for the cumulative Parisian ruin time $\tau_{ r }(u)$, %displayed in Theorem  \ref{Thm3}, 
The proof of Theorem \ref{Thm3} will be done by  first deriving the corresponding results for the cumulative Parisian ruin problem %time $\tau_{\widetilde{r}}(v)$
of the auxiliary risk model   \eqref{eq:Unew}, as $v\to\IF$, % (applying Theorem \ref{Thm1}),  
and then   using the equivalence described in 
\eqref{eq:tauuv}.

\def\walpha{\widetilde{\vk\alpha}}
\def\wmu{\widetilde{\vk \mu} }

%and thus analysis of the original problem with $u\to\IF$ is equivalent to the analysis of a new problem with $v\to\IF$
%thus finding the asymptotics of the cumulative Parisian ruin probability and the asymptotic distribution of the conditional P as $u\to\IF$ is equivalent to deriving the asymptotics of
%\BQNY
%\pk{\tau_{\widetilde{r}}(v)<\IF},\ \ \ v\to\IF,
%\EQNY
In order to apply Theorem \ref{Thm1} to the auxiliary risk model   \eqref{eq:Unew}, a crucial step is to find the minimiser of the  $g$-function   given by
\BQNY
g(t)=\frac{1}{t}\inf_{\vk{x}\ge \walpha +  \wmu  t} \vk{x}^\top\Sigma^{-1}\vk{x},   %\text{with}\ \%,\ \ \ \vk{b}_t=(1,b_t)^\top.
 \EQNY
for which we must first solve the quadratic programming problem $P_\Sigma( \walpha +\wmu t )$ involved. 
% for $d=2$.  
%Further, recall  in our notation  $I(t)$ denotes the  essential index  set of the quadratic programming problem $P_\Sigma( \vk\alpha +\vk \mu t )$. If $I(t)^c\neq \emptyset$  we define
%\BQN\label{eq:KK}
%K(t)=\{j\in I(t)^c:  \Sigma_{jI(t)}\Sigma_{I(t)I(t)}^{-1}(\vk\alpha +\vk \mu t)_{I(t)}=(\vk\alpha +\vk \mu t)_j\}.
%\EQN
To this end, we adopt a direct approach, % by using \nelem{AL}, 
which is different from that in \cite{DHJT18}. It follows from \nelem{AL} that the $g$-function has different expressions on different sets of $t$ defined below:
\begin{itemize}
\item[(S1).] On the set $E_1=\{t\ge0: \rho(\alpha +\mu  t)\ge ( 1+  t)\}$,\ \   %\inf_{\vk{v}\ge\vk{\mu} t+\vk 1 } \vk{v}^\top\Sigma^{-1}\vk{v}
$g(t)=g_2(t):=g_{\{2\}}(t) =\frac{1}{t}(\alpha +\mu  t)^2$;
\item [(S2).] On the set $E_2=\{t\ge0: \rho( 1+  t)\ge (\alpha +\mu  t)\}$,\ \   %\inf_{\vk{v}\ge\vk{\mu} t+\vk 1 } \vk{v}^\top\Sigma^{-1}\vk{v}
$g(t)=g_1(t):=g_{\{1\}}(t)= \frac{1}{t}(1+  t)^2$;
\item [(S3).] On the set $E_3 = [0,\IF) \setminus (E_1\cup E_2)$,\ \  $g(t)= g_0(t),$
\end{itemize}
where 
\BQNY
g_0(t):=g_{\{1,2\}}(t)&=&\frac{1}{t}(\walpha +\wmu t  )^\top \Sigma^{-1} (\walpha +\wmu  t ) \\
&=&
\frac{1+\alpha^2-2\alpha\rho}{1-\rho^2} \frac{1}{t} +\frac{2(1+\alpha\mu-\rho\alpha -\rho\mu)}{1-\rho^2} +\frac{1+\mu^2 -2\rho\mu}{1-\rho^2} t.
\EQNY
Moreover, it is easy to see that the unique minimiser of $g_0(t), t\ge0$ is $t_0^{(0)}=\sqrt{\frac{1+\alpha^2-2\alpha  \rho }{1+\mu ^2-2\mu \rho} }$, the unique minimiser of $g_1(t), t\ge0$ is   $t_0^{(1)}=1$, and the unique minimiser of $g_2(t), t\ge0$ is $t_0^{(2)}= {\alpha}/{\mu}$. Note that all functions $g_i(t), t\ge0$, $i=0,1,2,$ are decreasing to the left of their own minimiser and then increasing to infinity.

In order to find the global minimiser of the $g$-function and the exact form of the sets $E_i, i=1,2,3,$ we shall discuss the following four cases, separately.
\begin{itemize}
\item[(1).] $\mu<1$ and $\alpha<1$,
\quad (2). $\mu\ge 1$ and $\alpha\ge1$,
\item[(3).] $\mu<1$ and $\alpha\ge1$, \quad
 (4). $\mu\ge 1$ and $\alpha<1$.
\end{itemize}

\subsection{ Case (1) $\mu<1$ and $\alpha<1$ }  Clearly, we have $E_1=\emptyset$ for any $\rho\in (-1,1)$. To analyse $E_2$ we distinguish the following three sub-cases:
\BQNY
 (1.1). \ \alpha<\mu,\ \ \ (1.2).\ \alpha>\mu, \ \ \ (1.3).\ \alpha=\mu.
\EQNY
\subsubsection{\underline{Case (1.1) $\alpha<\mu$}} In this case, we have

\begin{eqnarray*}
E_2=  
\left\{
  \begin{array}{ll}
    \emptyset, & \hbox{if \ $-1< \rho  \le \alpha$;} \\
   \{t\ge 0: t\le w\}, & \hbox{if \ $\alpha  < \rho<\mu$ ;} \\
   {[0,\IF)} , & \hbox{if \ $\mu  \le\rho<1$,}
  \end{array}
\right.\ \ \ \ \text{with}\ \ w=\frac{\rho-\alpha}{\mu-\rho}.
\end{eqnarray*}
This combined with the fact that $E_1=\emptyset$ for any $\rho\in (-1,1)$ yields that, if  $-1< \rho  \le \alpha$ then $g(t)\equiv g_0(t), t\ge 0$ where the minimum is attained at the unique point $t_0^{(0)}$, if $\mu  \le\rho<1$ then  $g(t)\equiv g_1(t), t\ge 0$ where the minimum is attained at the unique point $\tzo$, and if $\alpha  < \rho<\mu$ then
\begin{eqnarray*} %\label{eq:g01}
g(t)=  
\left\{
  \begin{array}{ll}
    g_0(t), & \hbox{if \ $t>w$;} \\
   %\{t\ge 0: t\le w\}, & \hbox{if \ $\alpha  < \rho<\mu$ ;} \\
  g_1(t), & \hbox{if \ $t\le w$,}
  \end{array}
\right.  \ \ \ \text{with}\ g_0(w)=g_1(w),
\end{eqnarray*}
%with $g_0(w)=g_1(w)$ which is easy to check. 
%Clearly,  if $-1< \rho  \le \alpha$ then the the unique maximiser of $g(t), t\ge 0$ is $t_0^{(0)}$ and if $\mu  \le\rho<1$ the unique maximiser of $g(t), t\ge 0$ is $t_0^{(1)}$. 
%where minimiser is not clear so far.
%Next we focus on the case where $\alpha  < \rho<\mu$.  
%According to \eqref{eq:g01}, we have %in this case
and thus
$$
\inf_{t\ge0} g(t)=\min \LT(\inf_{t\le w} g_1(t), \ \inf_{t>w} g_0(t)\RT).
$$
In order to derive $\inf_{t\le w} g_1(t)$ and  $\inf_{t>w} g_0(t)$ we need to check if $\tzo<w$ and if $\tzz>w$. We can show that
 \BQN \label{eq:ar2}
 \tzo<w\ \ \Leftrightarrow\ \ \rho>\frac{\alpha+\mu}{2},\ \ \ \tzz>w \ \ \Leftrightarrow\ \ \rho<\frac{\alpha+\mu}{2},\ \ \ \tzz=\tzo=w \ \ \Leftrightarrow\ \ \rho=\frac{\alpha+\mu}{2}.
 \EQN 
Note that  $\alpha<\frac{\alpha+\mu}{2}<\mu$. Thus, we have  for $\alpha  < \rho<\mu$
\BQNY
\inf_{t\ge0} g(t)=\left\{
  \begin{array}{ll}
    g_0(\tzz), & \hbox{if \ $\alpha<\rho<\frac{\alpha+\mu}{2}$;} \\
     g_0(\tzz) = g_1(\tzo), & \hbox{if \ $\rho=\frac{\alpha+\mu}{2}$ ;} \\
  g_1(\tzo), & \hbox{if \ $\frac{\alpha+\mu}{2}<\rho<\mu$,}
  \end{array}
\right. 
\EQNY
and in each of the above three cases the minimiser of the function $g(t), t\ge0$ is unique. Using the notation in Theorem \ref{Thm1}, the above findings for Case (1.1) $\alpha<\mu$  are summarized in the following lemma:
 
\BEL \label{lem:two}
(1).  If $-1<\rho<\frac{\alpha+\mu}{2}$, then
\BQNY
t_0=\tzz,\ \ I= \{1,2\},\ \ K=\emptyset, \ \ \gt =g_0(\tzz) ,\ \   \ggt=g_0''(\tzz).
\EQNY

(2). If $ \rho=\frac{\alpha+\mu}{2}$, then
\BQNY
t_0=\tzz=\tzo=w,\ \ I= \{1\},\ \ K=\{2\},\ \  \gt=g_0(\tzz)=g_1(\tzo)=4,\ \   \ggt=  g_1''(\tzo)=2.
\EQNY

(3).  If $ \frac{\alpha+\mu}{2}<\rho<1$, then
\BQNY
t_0=\tzo,\ \ I= \{1\},\ \ K=\emptyset,\ \  \gt =g_1(\tzo)=4,\ \   \ggt=  g_1''(\tzo)=2.
\EQNY

\EEL

\subsubsection{\underline{Case (1.2) $\alpha>\mu$}} In this case, we have  (recall $ w=\frac{ \alpha-\rho}{\rho-\mu }$)
\begin{eqnarray*}
E_2=  
\left\{
  \begin{array}{ll}
    \emptyset, & \hbox{if \ $-1< \rho  \le \mu$;} \\
   \{t\ge 0: t\ge w\}, & \hbox{if \ $\mu  < \rho<\alpha$ ;} \\
   {[0,\IF)} , & \hbox{if \ $\alpha  \le\rho<1$.}
  \end{array}
\right. %\ \ \ \ \text{with}\ \ w=\frac{\rho-\alpha}{\mu-\rho}.
\end{eqnarray*}
%This combined with the fact that $E_1=\emptyset$ for any $\rho\in (-1,1)$ yields 
%Note that if $-1< \rho  \le \mu$ or $\alpha  \le\rho<1$
%Similarly as  in Case (1.1), 
Thus, we have that, if  $-1< \rho  \le \mu$ then $g(t)\equiv g_0(t), t\ge 0$ where the minimum is attained at the unique point $t_0^{(0)}$, if $\alpha  \le\rho<1$ then  $g(t)\equiv g_1(t), t\ge 0$ where the minimum is attained at the unique point $\tzo$, and if $\mu  < \rho<\alpha$ then
\begin{eqnarray*} %\label{eq:g01}
g(t)=  
\left\{
  \begin{array}{ll}
    g_0(t), & \hbox{if \ $t< w$;} \\
   %\{t\ge 0: t\le w\}, & \hbox{if \ $\alpha  < \rho<\mu$ ;} \\
  g_1(t), & \hbox{if \ $t\ge w$.}
  \end{array}
\right. % \ \ \ \text{with}\ g_0(w)=g_1(w),
\end{eqnarray*}
%According to \eqref{eq:g01}, we have %in this case
%and thus
%$$
%\inf_{t\ge0} g(t)=\min \LT(\inf_{t\ge w} g_1(t), \ \inf_{t<w} g_0(t)\RT).
%$$
Similarly as  in Case (1.1), we have  for $\mu  < \rho<\alpha$
\BQNY
\inf_{t\ge0} g(t)=\left\{
  \begin{array}{ll}
    g_0(\tzz), & \hbox{if \ $\mu<\rho<\frac{\alpha+\mu}{2}$;} \\
     g_0(\tzz) = g_1(\tzo), & \hbox{if \ $\rho=\frac{\alpha+\mu}{2}$ ;} \\
  g_1(\tzo), & \hbox{if \ $\frac{\alpha+\mu}{2}<\rho<\alpha$,}
  \end{array}
\right. 
\EQNY
and in each of the above three cases the minimiser of the function $g(t), t\ge0$ is unique.
Summarizing  the above we conclude that the results in \nelem{lem:two} still hold for Case (1.2) $\alpha>\mu$.

\subsubsection{\underline{Case (1.3) $\alpha=\mu$}} In this case, we have % (recall $ w=\frac{ \alpha-\rho}{\rho-\mu }$)
\begin{eqnarray*}
E_2=  
\left\{
  \begin{array}{ll}
    \emptyset, & \hbox{if \ $-1< \rho  <\mu$;} \\
  % \{t\ge 0: t\ge w\}, & \hbox{if \ $\mu  < \rho<\alpha$ ;} \\
   {[0,\IF)} , & \hbox{if \ $\mu  \le\rho<1$.}
  \end{array}
\right. %\ \ \ \ \text{with}\ \ w=\frac{\rho-\alpha}{\mu-\rho}.
\end{eqnarray*}
Using similar arguments as in Case (1.1) and noting that $\frac{\alpha+\mu}{2}=\alpha=\mu$, we could prove that 
the results in \nelem{lem:two} still hold for Case (1.3) $\alpha=\mu$.
\medskip

Consequently, we conclude that   \nelem{lem:two}  holds  for Case (1) $\mu<1$ and $\alpha<1$.

\subsection{Case (2) $\mu\ge 1$ and $\alpha\ge1$}  Clearly, we have $E_2=\emptyset$ for any $\rho\in (-1,1)$. To analyse $E_1$ we distinguish the following three sub-cases:
\BQNY
 (2.1). \ \alpha<\mu,\ \ \ (2.2).\ \alpha>\mu, \ \ \ (2.3).\ \alpha=\mu.
\EQNY
\subsubsection{\underline{Case (2.1) $\alpha<\mu$}} In this case, we have  
\begin{eqnarray*}
E_1=  
\left\{
  \begin{array}{ll}
    \emptyset, & \hbox{if \ $-1< \rho  \le 1/\mu$;} \\
   \{t\ge 0: t\ge Q\}, & \hbox{if \ $1/\mu  < \rho<1/\alpha$ ;} \\
   {[0,\IF)} , & \hbox{if \ $1/\alpha  \le\rho<1$,}
  \end{array}
\right.\ \ \ \ \text{with}\ \ Q=\frac{1-\rho\alpha}{\mu\rho-1}.
\end{eqnarray*}
This combined with the fact that $E_2=\emptyset$ for any $\rho\in (-1,1)$ yields that, if  $-1< \rho  \le 1/\mu$ then $g(t)\equiv g_0(t), t\ge 0$ where the minimum is attained at the unique point $t_0^{(0)}$, if $1/\alpha \le\rho<1$ then  $g(t)\equiv g_2(t), t\ge 0$ where the minimum is attained at the unique point $\tzt$, and if $1/\mu  < \rho<1/\alpha$ then
\begin{eqnarray*} %\label{eq:g01}
g(t)=  
\left\{
  \begin{array}{ll}
    g_0(t), & \hbox{if \ $t<Q$;} \\
   %\{t\ge 0: t\le w\}, & \hbox{if \ $\alpha  < \rho<\mu$ ;} \\
  g_2(t), & \hbox{if \ $t\ge Q$,}
  \end{array}
\right.  \ \ \ \text{with}\ g_0(Q)=g_2(Q),
\end{eqnarray*}
%with $g_0(w)=g_1(w)$ which is easy to check. 
%Clearly,  if $-1< \rho  \le \alpha$ then the the unique maximiser of $g(t), t\ge 0$ is $t_0^{(0)}$ and if $\mu  \le\rho<1$ the unique maximiser of $g(t), t\ge 0$ is $t_0^{(1)}$. 
%where minimiser is not clear so far.
%Next we focus on the case where $\alpha  < \rho<\mu$.  
%According to \eqref{eq:g01}, we have %in this case
and thus
$$
\inf_{t\ge0} g(t)=\min \LT(\inf_{t\ge Q} g_2(t), \ \inf_{t<Q} g_0(t)\RT).
$$
In order to derive $\inf_{t\ge Q} g_2(t)$ and  $\inf_{t<Q} g_0(t)$ we need to check if $\tzt>Q$ and if $\tzz<Q$. We can show that
 \BQN \label{eq:arQ}
 \tzt>Q\ \ \Leftrightarrow\ \ \rho>\frac{\alpha+\mu}{2\alpha\mu},\ \ \ \tzz<Q \ \ \Leftrightarrow\ \ \rho<\frac{\alpha+\mu}{2\alpha\mu},\ \ \ \tzz=\tzt=Q \ \ \Leftrightarrow\ \ \rho=\frac{\alpha+\mu}{2\alpha\mu}.
 \EQN 
Note that  $1/\mu<\frac{\alpha+\mu}{2\alpha\mu}<1/\alpha$. Thus, we have  for $1/\mu< \rho<1/\alpha$
\BQNY
\inf_{t\ge0} g(t)=\left\{
  \begin{array}{ll}
    g_0(\tzz), & \hbox{if \ $1/\mu<\rho<\frac{\alpha+\mu}{2\alpha\mu}$;} \\
     g_0(\tzz) = g_2(\tzt), & \hbox{if \ $\rho=\frac{\alpha+\mu}{2\alpha\mu}$ ;} \\
  g_2(\tzt), & \hbox{if \ $\frac{\alpha+\mu}{2\alpha\mu}<\rho<1/\alpha$,}
  \end{array}
\right. 
\EQNY
and in each of the above three cases the minimiser of the function $g(t), t\ge0$ is unique. Using the notation in Theorem \ref{Thm1}, the above findings for Case (2.1)  $\alpha<\mu$  are summarized in the following lemma:
 
\BEL \label{lem:two-2}
(1).  If $-1<\rho<\frac{\alpha+\mu}{2\alpha\mu}$, then
\BQNY
t_0=\tzz,\ \ I= \{1,2\},\ \ K=\emptyset, \ \ \gt =g_0(\tzz) ,\ \   \ggt=g_0''(\tzz).
\EQNY

(2). If $ \rho=\frac{\alpha+\mu}{2\alpha\mu}$, then
\BQNY
t_0=\tzz=\tzt=Q,\ \ I= \{2\},\ \ K=\{1\},\ \  \gt=g_0(\tzz)=g_2(\tzt)=4\alpha\mu,\ \   \ggt=  g_2''(\tzt)=2\alpha^{-1}\mu^3.
\EQNY

(3).  If $ \frac{\alpha+\mu}{2\alpha\mu}<\rho<1$, then
\BQNY
t_0=\tzt,\ \ I= \{2\},\ \ K=\emptyset,\ \  \gt =g_2(\tzt)=4\alpha\mu,\ \   \ggt=  g_2''(\tzt)=2\alpha^{-1}\mu^3.
\EQNY

\EEL

Case (2.1) and Case (2.2) can be analysed similarly, and %similarly as in Section 6.1 
we can conclude that \nelem{lem:two-2} holds for Case (2) $\mu\ge 1$ and $\alpha\ge1$.

\subsection{Case (3) $\mu<1$ and $\alpha\ge1$.}  
%Clearly, we have $E_2=\emptyset$ for any $\rho\in (-1,1)$.
%First note that for the constants $w, Q$ given in Case (1) and Case (2) we can check that $Q<w$ in this case. 
To analyse $E_1, E_2$ we distinguish the following three sub-cases:
\BQNY
 (3.1). \ \mu<1/\alpha,\ \ \ (3.2).\ \mu>1/\alpha, \ \ \ (3.3).\ \mu=1/\alpha.
\EQNY
\subsubsection{\underline{Case (3.1) $\mu<1/\alpha$}} In this case, we have  
\begin{eqnarray}\label{eq:E12}
E_1=  
\left\{
  \begin{array}{ll}
    \emptyset, & \hbox{if \ $-1< \rho  \le 1/\alpha$;} \\
   \{t\ge 0: t\le Q\}, & \hbox{if \ $1/\alpha  < \rho<1$,}  %\\
  % {[0,\IF)} , & \hbox{if \ $\alpha  \le\rho<1$,}
  \end{array}
\right.\ \  E_2=  
\left\{
  \begin{array}{ll}
    \emptyset, & \hbox{if \ $-1< \rho  \le \mu$;} \\
  % \{t\ge 0: t\ge Q\}, & \hbox{if \ $1/\mu  < \rho<\alpha$ ;} \\
  \{t\ge 0: t\ge w\}, & \hbox{if \ $\mu  <\rho<1$.}
  \end{array}
\right. 
\end{eqnarray}
This  implies that, if  $-1< \rho  \le \mu$ then $g(t)\equiv g_0(t), t\ge 0,$ where the minimum is attained at the unique point $t_0^{(0)}$,  if $\mu  < \rho \le 1/\alpha$ then
\begin{eqnarray*} %\label{eq:g01}
g(t)=  
\left\{
  \begin{array}{ll}
    g_0(t), & \hbox{if \ $t<w$;} \\
   %\{t\ge 0: t\le w\}, & \hbox{if \ $\alpha  < \rho<\mu$ ;} \\
  g_1(t), & \hbox{if \ $t\ge w$,}
  \end{array}
\right. % \ \ \ \text{with}\ g_0(Q)=g_2(Q),
\end{eqnarray*}
implying that 
$$
\inf_{t\ge0} g(t)=\min \LT(\inf_{t\ge w} g_1(t), \ \inf_{t<w} g_0(t)\RT),
$$
and if $1/\alpha <\rho<1$ then  
\begin{eqnarray*} %\label{eq:g01}
g(t)=  
\left\{
  \begin{array}{ll}
    g_2(t), & \hbox{if \ $t\le Q$;} \\
  g_0(t), & \hbox{if \ $Q  < t<w$ ;} \\
  g_1(t), & \hbox{if \ $t\ge w$.}
  \end{array}
\right. % \ \ \ \text{with}\ g_0(Q)=g_2(Q),
\end{eqnarray*}
implying that
\BQN\label{eq:g012-0}
\inf_{t\ge0} g(t)=\min \LT(\inf_{t\le Q} g_2(t),   \ \inf_{Q<t<w} g_0(t),\ \inf_{t\ge w} g_1(t)\RT).
\EQN
Note that  $Q<w$ for any $\mu<\rho<1$. 
%In order to obtain the values of $\inf_{t\ge Q} g_2(t)$ and  $\inf_{t<Q} g_0(t)$ we have to check if $\tzt>Q$ and if $\tzz<Q$. We can show that
Similarly as \eqref{eq:ar2} and \eqref{eq:arQ} we have, for any  $\mu  < \rho<1 (\le \alpha)$ 
 \BQN  \label{eq:tzto}
 \tzo<w\ \ \Leftrightarrow\ \ \rho<\frac{\alpha+\mu}{2},\ \ \ \tzz>w \ \ \Leftrightarrow\ \ \rho>\frac{\alpha+\mu}{2},\ \ \ \tzz=\tzo=w \ \ \Leftrightarrow\ \ \rho=\frac{\alpha+\mu}{2},
\EQN 
and, for any $1/\alpha <\rho<1$
\BQN\label{eq:tztz}
\tzt<Q\ \ \Leftrightarrow\ \ \rho>\frac{\alpha+\mu}{2\alpha\mu},\ \ \ \tzz>Q \ \ \Leftrightarrow\ \ \rho<\frac{\alpha+\mu}{2\alpha\mu},\ \ \ \tzz=\tzt=Q \ \ \Leftrightarrow\ \ \rho=\frac{\alpha+\mu}{2\alpha\mu}.
 \EQN 
Furthermore, it follows that  
$\frac{\alpha+\mu}{2\alpha\mu}> \frac{1}{2}(1/\mu+\mu)> 1$ in the considered case $\mu<1/\alpha\le 1$, then from \eqref{eq:tztz} we conclude that $\tzt>Q,  \tzz>Q$ hold  for  any $1/\alpha <\rho<1$,  which helps to further simplify \eqref{eq:g012-0} as follows,
 \BQN
\inf_{t\ge0} g(t)=\min \LT(\inf_{t\le Q} g_2(t),   \ \inf_{Q<t<w} g_0(t),\ \inf_{t\ge w} g_1(t)\RT)=\min \LT( \inf_{Q<t<w} g_0(t),\ \inf_{t\ge w} g_1(t)\RT).
\EQN
After some simple calculations as before, we can show that  for $\mu< \rho<1$
\BQNY
\inf_{t\ge0} g(t)=\left\{
  \begin{array}{ll}
    g_0(\tzz), & \hbox{if \ $ \mu<\rho<\frac{\alpha+\mu}{2}$;} \\
     g_0(\tzz) = g_1(\tzo), & \hbox{if \ $\rho=\frac{\alpha+\mu}{2}$ ;} \\
  g_1(\tzo), & \hbox{if \ $\frac{\alpha+\mu}{2}<\rho<1$,}
  \end{array}
\right. 
\EQNY
and in each of the above three cases the minimiser of the function $g(t), t\ge0$ is unique. Summarizing the above findings we can conclude that \nelem{lem:two} holds for Case (3.1).

\subsubsection{\underline{Case (3.2) $\mu>1/\alpha$}} In this case, we have  \eqref{eq:E12} still holds.
It follows  that, if  $-1< \rho  \le 1/\alpha$ then $g(t)\equiv g_0(t), t\ge 0,$ where the minimum is attained at the unique point $t_0^{(0)}$,  if $1/\alpha < \rho \le \mu$ then
\begin{eqnarray*} %\label{eq:g01}
g(t)=  
\left\{
  \begin{array}{ll}
    g_0(t), & \hbox{if \ $t>Q$;} \\
   %\{t\ge 0: t\le w\}, & \hbox{if \ $\alpha  < \rho<\mu$ ;} \\
  g_2(t), & \hbox{if \ $t\le Q$,}
  \end{array}
\right. % \ \ \ \text{with}\ g_0(Q)=g_2(Q),
\end{eqnarray*}
implying that 
$$
\inf_{t\ge0} g(t)=\min \LT(\inf_{t\le Q} g_2(t), \ \inf_{t>Q} g_0(t)\RT),
$$
and if $\mu <\rho<1$ then  
\begin{eqnarray*} %\label{eq:g01}
g(t)=  
\left\{
  \begin{array}{ll}
    g_2(t), & \hbox{if \ $t\le Q$;} \\
  g_0(t), & \hbox{if \ $Q  < t<w$ ;} \\
  g_1(t), & \hbox{if \ $t\ge w$.}
  \end{array}
\right. % \ \ \ \text{with}\ g_0(Q)=g_2(Q),
\end{eqnarray*}
implying that
\BQN\label{eq:g012-1}
\inf_{t\ge0} g(t)=\min \LT(\inf_{t\le Q} g_2(t),   \ \inf_{Q<t<w} g_0(t),\ \inf_{t\ge w} g_1(t)\RT).
\EQN
%In order to obtain the values of $\inf_{t\ge Q} g_2(t)$ and  $\inf_{t<Q} g_0(t)$ we have to check if $\tzt>Q$ and if $\tzz<Q$. We can show that
Note that in this case both \eqref{eq:tzto} and \eqref{eq:tztz} are still valid for the corresponding values of $\rho$ mentioned therein.
Furthermore, it follows that  
$\frac{\alpha+\mu}{2}\ge \sqrt{\alpha\mu}> 1$ in the considered case $\mu>1/\alpha$, then from \eqref{eq:tzto} we conclude that $\tzo<w,  \tzz<w$ hold  for  any $\mu <\rho<1$,  which helps to further simplify \eqref{eq:g012-1} as follows,
 \BQN
\inf_{t\ge0} g(t)=\min \LT(\inf_{t\le Q} g_2(t),   \ \inf_{Q<t<w} g_0(t),\ \inf_{t\ge w} g_1(t)\RT)=\min \LT( \inf_{t\le Q} g_2(t),    \inf_{Q<t<w} g_0(t)\RT).
\EQN
Thus,  we can show that for $1/\alpha< \rho<1$
\BQNY
\inf_{t\ge0} g(t)=\left\{
  \begin{array}{ll}
    g_0(\tzz), & \hbox{if \ $ 1/\alpha<\rho<\frac{\alpha+\mu}{2\alpha\mu}$;} \\
     g_0(\tzz) = g_2(\tzt), & \hbox{if \ $\rho=\frac{\alpha+\mu}{2\alpha\mu}$ ;} \\
  g_2(\tzt), & \hbox{if \ $\frac{\alpha+\mu}{2\alpha\mu}<\rho<1$,}
  \end{array}
\right. 
\EQNY
and in each of the above three cases the minimiser of the function $g(t), t\ge0$ is unique. Summarizing the above findings we can conclude that \nelem{lem:two-2} holds for Case (3.2).

\subsubsection{\underline{Case (3.3) $\mu=1/\alpha$}} In this case,  
we have  that \eqref{eq:E12} still holds. As now $\frac{\alpha+\mu}{2}=\frac{\alpha+\mu}{2\alpha\mu}\ge 1$, we obtain that (i) in \nelem{lem:two} (the same to (i) in \nelem{lem:two-2}) is valid for any $-1<\rho<1$.
%This   yields that, if  $-1< \rho  \le 1/\alpha$ then $g(t)\equiv g_0(t), t\ge 0$ where the minimum is attained at the unique point $t_0^{(0)}$
\medskip

Consequently, we can conclude that for Case (3) $\mu< 1$ and $\alpha\ge1$, if further $\mu\le 1/\alpha$ then \nelem{lem:two} holds, and if further $\mu> 1/\alpha$ then \nelem{lem:two-2} holds.
 \subsection{Case (4) $\mu\ge 1$ and $\alpha<1$.}
 To analyse $E_1, E_2$ we distinguish the following three sub-cases:
\BQNY
 (4.1). \ 1/\mu<\alpha,\ \ \ (4.2).\ 1/\mu>\alpha, \ \ \ (4.3).\ 1/\mu=\alpha.
\EQNY
\subsubsection{\underline{Case (4.1) $1/\mu<\alpha$}} In this case, we have  
\begin{eqnarray}\label{eq:E12-2}
E_1=  
\left\{
  \begin{array}{ll}
    \emptyset, & \hbox{if \ $-1< \rho  \le 1/\mu$;} \\
   \{t\ge 0: t\ge Q\}, & \hbox{if \ $1/\mu  < \rho<1$,}  %\\
  % {[0,\IF)} , & \hbox{if \ $\alpha  \le\rho<1$,}
  \end{array}
\right.\ \  E_2=  
\left\{
  \begin{array}{ll}
    \emptyset, & \hbox{if \ $-1< \rho  \le \alpha$;} \\
  % \{t\ge 0: t\ge Q\}, & \hbox{if \ $1/\mu  < \rho<\alpha$ ;} \\
  \{t\ge 0: t\le w\}, & \hbox{if \ $\alpha  <\rho<1$.}
  \end{array}
\right. 
\end{eqnarray}
This  implies that, if  $-1< \rho  \le 1/\mu$ then $g(t)\equiv g_0(t), t\ge 0$, where the minimum is attained at the unique point $t_0^{(0)}$,  if $1/\mu  < \rho \le  \alpha$ then
\begin{eqnarray*} %\label{eq:g01}
g(t)=  
\left\{
  \begin{array}{ll}
    g_0(t), & \hbox{if \ $t<Q$;} \\
   %\{t\ge 0: t\le w\}, & \hbox{if \ $\alpha  < \rho<\mu$ ;} \\
  g_2(t), & \hbox{if \ $t\ge Q$,}
  \end{array}
\right. % \ \ \ \text{with}\ g_0(Q)=g_2(Q),
\end{eqnarray*}
implying that 
$$
\inf_{t\ge0} g(t)=\min \LT(\inf_{t\ge Q} g_2(t), \ \inf_{t<Q} g_0(t)\RT),
$$
and if $\alpha <\rho<1$ then  
\begin{eqnarray*} %\label{eq:g01}
g(t)=  
\left\{
  \begin{array}{ll}
    g_2(t), & \hbox{if \ $t\ge Q$;} \\
  g_0(t), & \hbox{if \ $w  < t<Q$ ;} \\
  g_1(t), & \hbox{if \ $t\le w$.}
  \end{array}
\right. % \ \ \ \text{with}\ g_0(Q)=g_2(Q),
\end{eqnarray*}
implying that
\BQN\label{eq:g012}
\inf_{t\ge0} g(t)=\min \LT(\inf_{t\ge Q} g_2(t),   \ \inf_{w<t<Q} g_0(t),\ \inf_{t\le w} g_1(t)\RT).
\EQN
Note that  $Q>w$ for any $1/\mu<\rho<1$. 
%In order to obtain the values of $\inf_{t\ge Q} g_2(t)$ and  $\inf_{t<Q} g_0(t)$ we have to check if $\tzt>Q$ and if $\tzz<Q$. We can show that
Similarly as in Case (3.1) we have  for $1/\mu< \rho<1$
\BQNY
\inf_{t\ge0} g(t)=\left\{
  \begin{array}{ll}
    g_0(\tzz), & \hbox{if \ $ 1/\mu<\rho<\frac{\alpha+\mu}{2\alpha\mu}$;} \\
     g_0(\tzz) = g_2(\tzt), & \hbox{if \ $\rho=\frac{\alpha+\mu}{2\alpha\mu}$ ;} \\
  g_2(\tzt), & \hbox{if \ $\frac{\alpha+\mu}{2\alpha\mu}<\rho<1$,}
  \end{array}
\right. 
\EQNY
and in each of the above three cases the minimiser of the function $g(t), t\ge0$ is unique. Summarizing the above findings we can conclude that \nelem{lem:two-2} holds for Case (4.1).

\subsubsection{\underline{Case (4.2) $1/\mu>\alpha$}} In this case, we have  \eqref{eq:E12-2} still holds.
It follows  that, if  $-1< \rho  \le  \alpha$ then $g(t)\equiv g_0(t), t\ge 0$, where the minimum is attained at the unique point $t_0^{(0)}$,  if $
\alpha < \rho \le 1/\mu$ then
\begin{eqnarray*} %\label{eq:g01}
g(t)=  
\left\{
  \begin{array}{ll}
    g_0(t), & \hbox{if \ $t>w$;} \\
   %\{t\ge 0: t\le w\}, & \hbox{if \ $\alpha  < \rho<\mu$ ;} \\
  g_1(t), & \hbox{if \ $t\le w$,}
  \end{array}
\right. % \ \ \ \text{with}\ g_0(Q)=g_2(Q),
\end{eqnarray*}
implying that 
$$
\inf_{t\ge0} g(t)=\min \LT(\inf_{t\le w} g_1(t), \ \inf_{t>w} g_0(t)\RT),
$$
and if $1/\mu <\rho<1$ then  
\begin{eqnarray*} %\label{eq:g01}
g(t)=  
\left\{
  \begin{array}{ll}
    g_2(t), & \hbox{if \ $t\ge Q$;} \\
  g_0(t), & \hbox{if \ $w  < t<Q$ ;} \\
  g_1(t), & \hbox{if \ $t\le w$.}
  \end{array}
\right. % \ \ \ \text{with}\ g_0(Q)=g_2(Q),
\end{eqnarray*}
implying that
\BQN\label{eq:g012-2}
\inf_{t\ge0} g(t)=\min \LT(\inf_{t\ge Q} g_2(t),   \ \inf_{w<t<Q} g_0(t),\ \inf_{t\le w} g_1(t)\RT).
\EQN
%In order to obtain the values of $\inf_{t\ge Q} g_2(t)$ and  $\inf_{t<Q} g_0(t)$ we have to check if $\tzt>Q$ and if $\tzz<Q$. We can show that
Similarly as before,  we can show that for $\alpha< \rho<1$
\BQNY
\inf_{t\ge0} g(t)=\left\{
  \begin{array}{ll}
    g_0(\tzz), & \hbox{if \ $ \alpha<\rho<\frac{\alpha+\mu}{2}$;} \\
     g_0(\tzz) = g_1(\tzo), & \hbox{if \ $\rho=\frac{\alpha+\mu}{2}$ ;} \\
  g_1(\tzo), & \hbox{if \ $\frac{\alpha+\mu}{2}<\rho<1$,}
  \end{array}
\right. 
\EQNY
and in each of the above three case the minimiser of the function $g(t), t\ge0$ is unique. Summarizing the above findings we can conclude that \nelem{lem:two} holds for Case (4.2).

\subsubsection{\underline{Case (4.3) $1/\mu=\alpha$}} In this case, we have  that
\eqref{eq:E12-2} still holds.
%We can follow the same arguments as in Case (3.1). 
As now $\frac{\alpha+\mu}{2}=\frac{\alpha+\mu}{2\alpha\mu}\ge 1$, we obtain that (i) in \nelem{lem:two-2} (the same to (i) in \nelem{lem:two}) is valid for any $-1<\rho<1$.
%This   yields that, if  $-1< \rho  \le 1/\alpha$ then $g(t)\equiv g_0(t), t\ge 0$ where the minimum is attained at the unique point $t_0^{(0)}$

\medskip

Consequently, we can conclude that for Case (4) $\mu\ge 1$ and $\alpha<1$, if further $1/\mu< \alpha$ then \nelem{lem:two-2} holds, and if further $1/\mu\ge  \alpha$ then \nelem{lem:two} holds.

\medskip

The following corollary is a collection of the main findings in Sections 6.1, 6.2, 6.3, 6.4.

\BK \label{Corol} 
Consider the auxiliary  two-dimensional Brownian motion models described in \eqref{eq:Unew}.
%Recall the notation in \eqref{eq:uvr}. 
\begin{itemize}

\item[(i).] Suppose $\alpha$ and $\mu$ satisfy one of the following conditions:
\begin{itemize}
\item[(i.C1)] $\mu<1$ and $\alpha<1$,
\item[(i.C2)] $\mu<1$, $\alpha\ge 1$ and $\mu\le 1/\alpha$,
\item[(i.C3)]  $\mu\ge 1$, $\alpha<1$ and $\mu\le 1/\alpha$.
\end{itemize}
We have %, for  any $r\ge0,$  

\begin{itemize}
\item[(i.R1).] If $-1<\rho < \frac{\alpha+\mu}{2}$,  then %, as $u\to\IF$
\BQNY
t_0=\tzz,\ \ I= \{1,2\},\ \ K=\emptyset, \ \ \gt =g_0(\tzz) ,\ \   \ggt=g_0''(\tzz).
\EQNY
\item[(i.R2).] If $ \rho=\frac{\alpha+\mu}{2 }$,  then %, as $u\to\IF$
\BQNY
t_0=\tzz=\tzo,\ \ I= \{1\},\ \ K=\{2\},\ \  \gt=g_0(\tzz)=g_1(\tzo)=4,\ \   \ggt=  g_1''(\tzo)=2.
\EQNY

\item[(i.R3).]  If $\frac{\alpha+\mu}{2}<\rho<1$,  then %, as $u\to\IF$
 \BQNY
t_0=\tzo,\ \ I= \{1\},\ \ K=\emptyset,\ \  \gt =g_1(\tzo)=4,\ \   \ggt=  g_1''(\tzo)=2.
\EQNY
\end{itemize}

\item[(ii).]  Suppose $\alpha$ and $\mu$  satisfy one of the following conditions:
\begin{itemize}
\item[(ii.C1)] $\mu\ge 1$ and $\alpha\ge1$,
\item[(ii.C2)] $\mu<1$, $\alpha\ge 1$ and $\mu> 1/\alpha$,
\item[(ii.C3)]  $\mu\ge 1$, $\alpha<1$ and $\mu> 1/\alpha$.
\end{itemize}
We have 
\begin{itemize}
\item[(ii.R1).] If $-1<\rho < \frac{\alpha+\mu}{2\alpha\mu}$, then 
 \BQNY
t_0=\tzz,\ \ I= \{1,2\},\ \ K=\emptyset, \ \ \gt =g_0(\tzz) ,\ \   \ggt=g_0''(\tzz).
\EQNY

\item[(ii.R2).] If $ \rho=\frac{\alpha+\mu}{2 \alpha \mu}$,  then
\BQNY
t_0=\tzz=\tzt,\ \ I= \{2\},\ \ K=\{1\},\ \  \gt=g_0(\tzz)=g_2(\tzt)=4\alpha\mu,\ \   \ggt=  g_2''(\tzt)=2\alpha^{-1}\mu^3.
\EQNY

\item[(i.R3).]  If $\frac{\alpha+\mu}{2}<\rho<1$,  then 
 \BQNY
t_0=\tzt,\ \ I= \{2\},\ \ K=\emptyset,\ \  \gt =g_2(\tzt)=4\alpha\mu,\ \   \ggt=  g_2''(\tzt)=2\alpha^{-1}\mu^3.
\EQNY

\end{itemize}
\end{itemize}

\EK

%\medskip

Consequently, by using results in Corollary \ref{Corol} and applying Theorem \ref{Thm1} we can obtain asymptotic results for the cumulative Parisian ruin probability and the conditional cumulative Parisian ruin time of
 the auxiliary risk model  \eqref{eq:Unew}, as $v\to\IF$. Finally, Theorem \ref{Thm3} follows by directly using the equivalence described in 
\eqref{eq:tauuv}.

\medskip

{\bf Acknowledgement}: The author would like to thank the referees for their carefully reading and constructive suggestions which
significantly improved the manuscript.
The author would also like to thank % Hansj\"{o}rg Albrecher, 
 Krzysztof D\c{e}bicki, Enkelejd Hashorva  and Tomasz Rolski
for their stimulating discussions on the topic of this paper.

\bibliographystyle{siam} %{ieeetr} %

 \bibliography{vectProcEKEEKK}
\end{document}